\documentclass[11pt, a4paper, twoside]{amsart}
\usepackage{amsmath}
\usepackage{amssymb}
\usepackage{amsfonts}
\usepackage{a4}
\usepackage[latin1]{inputenc}

\numberwithin{equation}{section}

\newtheorem{theorem}{Theorem}[section]

\theoremstyle{definition} 

\theoremstyle{plain} 

\newtheorem{conclusion}[theorem]{Conclusion}

\newtheorem{claim}[theorem]{Claim}

\newtheorem*{maintheorem*}{Main Theorem} 
\newtheorem*{conjecture*}{Conjecture} 
\newtheorem{definition}[theorem]{Definition}

\theoremstyle{remark}  
\newtheorem*{remarks*}{Remarks}
\newtheorem*{remark*}{Remark}
\newtheorem*{claim*}{Claim}
\newcommand{\nc}{\newcommand}
\nc{\nothing}[1]{}

\nothing{   
\nc{\dom}{{\rm dom}}
\nc{\card}{{\rm card}}
\nc{\lh}{{\rm lh}}
\nc{\lgg}{{\rm lg}}
\nc{\rge}{\mbox{\rm range}}
\nc{\cf}{{\rm cf}}
\nc{\nex}{\mbox{\rm next}}
\nc{\uhr}{\restriction}
\nc{\supt}{{\rm supt}}
\nc{\supp}{{\rm supp}}
\nc{\Lim}{{\rm Lim}}
\nc{\Leb}{{\rm Leb}}
\nc{\modd}{{\rm mod}}
\nc{\RO}{{\rm RO}}
\nc{\prob}{{\rm Prob}}
}

\nc{\On}{{\rm On}}

\nc{\nco}{\DeclareMathOperator}

\nco{\order}{o} 
\nco{\ppower}{pp} 
\nco{\pcf}{pcf} 
\nco{\tcf}{tcf} 
\nco{\tlim}{tlim} 
\nco{\limtext}{lim} 
\nco{\prodt}{{\textstyle \prod}} 
\nco{\symdiff}{\triangle}
\nco{\dom}{dom}
\nco{\card}{card}
\nco{\lh}{lh}
\nco{\lgg}{lg}
\nco{\rge}{range}
\nco{\otp}{otp}
\nco{\trunk}{tr}
\nco{\cf}{cf}
\nco{\nex}{next}
\nc{\uhr}{\restriction}
\nco{\supt}{supt}
\nco{\supp}{supp}
\nco{\Lim}{Lim}
\nco{\Leb}{Leb}
\nco{\modd}{mod}
\nco{\invariant}{inv}
\nco{\id}{id}
\nco{\RO}{RO}

\nc{\potom}{\ensuremath{{\cal P}(\omega)}}
\nc{\potinf}{\ensuremath{[\omega]^\omega}}
\nc{\pfin}{\ensuremath{{\cal P}(\omega)/{\rm fin}}}

\nc{\potfin}{\ensuremath{[\omega]^{<\omega}}}
\nc{\inn}{\ensuremath{{\omega^{\uparrow \omega}}}}
\nc{\hoch}{^{<\omega}}
\nc{\hocho}{^{\omega}}
\nc{\tree}[1]{{[} #1 {]}_0}
\nc{\tre}[2]{ {#1}_{#2}}
   
\nc{\prooff}[1]{{\bf Proof} of #1:}
\nc{\proofend}{\makebox{} \hfill ${\bf \square}$ \\}
\nc{\proofendof}[1]{\makebox{} \hfill $\boldmath{\square}_{\rm #1}$ \\}
\nc{\beq}{\begin{eqnarray*}}
\nc{\eeq}{\end{eqnarray*}}
\nc{\bde}{\begin{list}}
\nc{\ede}{\end{list}}

\newenvironment{myrules}
{\begin{list}{}
{
 \setlength{\leftmargin}{0.5in}
 \setlength{\labelwidth}{1cm}
 \setlength{\labelsep}{0.2in}
 \setlength{\parsep}{0.5ex plus 0.2ex minus 0.1 ex}
 \setlength{\itemsep}{0.3ex plus 0.2 ex minus 0ex}
}}{\end{list}}

{\end{list}}

\newcounter{subalph}
{\end{list}}

\newcommand{\greek}[1]{\ifthenelse{\value{#1}=1}{\mbox{$\alpha$}}%
  {\ifthenelse{\value{#1}=2}{\mbox{$\beta$}}{%
   \ifthenelse{\value{#1}=3}{\mbox{$\gamma$}}{%
   \ifthenelse{\value{#1}=4}{\mbox{$\delta$}}{%
   \ifthenelse{\value{#1}=5}{\mbox{$\varepsilon$}}{%
   \ifthenelse{\value{#1}=6}{\mbox{$\zeta$}}{%
   \ifthenelse{\value{#1}=7}{\mbox{$\eta$}}{%
   \ifthenelse{\value{#1}=8}{\mbox{$\theta$}}{%
   \ifthenelse{\value{#1}=9}{\mbox{$\iota$}}{%
   \ifthenelse{\value{#1}=10}{\mbox{$\kappa$}}{%
   \ifthenelse{\value{#1}=11}{\mbox{$\lambda$}}{%
   \ifthenelse{\value{#1}=12}{\mbox{$\mu$}}{%
   \ifthenelse{\value{#1}=13}{\mbox{$\nu$}}{%
   \ifthenelse{\value{#1}=14}{\mbox{$\xi$}}{%
   \ifthenelse{\value{#1}=15}{\mbox{$\rm o$}}{%
   \ifthenelse{\value{#1}=16}{\mbox{$\pi$}}{%
   \ifthenelse{\value{#1}=17}{\mbox{$\varrho$}}{%
   \ifthenelse{\value{#1}=18}{\mbox{$\sigma$}}{%
   \ifthenelse{\value{#1}=19}{\mbox{$\tau$}}{%
   \ifthenelse{\value{#1}=20}{\mbox{$\upsilon$}}{%
   \ifthenelse{\value{#1}=21}{\mbox{$\varphi$}}{%
   \ifthenelse{\value{#1}=22}{\mbox{$\chi$}}{%
   \ifthenelse{\value{#1}=23}{\mbox{$\psi$}}{\mbox{$\omega$}%
  }}}}}}}}}}}}}}}}}}}}}}}}

\newcounter{subgreek}
{\end{list}}

\newcounter{subarabic}
{\end{list}}

\newcounter{subroman}
{\end{list}}

\newcount\skewfactor

\def\mathunderaccent#1#2 {\let\theaccent#1\skewfactor#2
\mathpalette\putaccentunder}
\def\putaccentunder#1#2{\oalign{$#1#2$\crcr\hidewidth
\vbox to.2ex{\hbox{$#1\skew\skewfactor\theaccent{}$}\vss}\hidewidth}}
\def\name{\mathunderaccent\tilde-3 }

\nc{\nname}{\name}

\nc{\even}{\ensuremath{\rm Even}}
\nc{\odd}{\ensuremath{\rm Odd}}

\nc{\al}{$\alpha$\  }
\nc{\om}{\omega}
\nc{\omm}{\ensuremath{\omega_1}}
\nc{\ep}{\varepsilon}
\nc{\tk}{\tilde{K}}
\nc{\concat}{\,\hat{} \,}   
\nc{\force}{\Vdash}
\nc{\fb}{f_{\bar{M}}}
\nc{\such}{\, : \,}   

\nc{\meager}{\ensuremath{{\cal M}}}
\nc{\lebesgue}{\ensuremath{{\cal N}}}
\nc{\nulll}{\ensuremath{{\cal N}}}
\nc{\ksigma}{\ensuremath{{\bf K}_\sigma}}
\nc{\ideal}{\ensuremath{{\cal I}}}
\nc{\ga}{\ensuremath{\frak a}}
\nc{\AAA}{{\cal A}}   
\nc{\gc}{\ensuremath{\frak c}}
\nc{\gs}{\ensuremath{\frak s}}
\nc{\gh}{\ensuremath{\frak h}}
\nc{\gd}{\ensuremath{\frak d}}
\nc{\gb}{\ensuremath{\frak b}}
\nc{\gro}{\ensuremath{\frak g}}
\nc{\gu}{\ensuremath{\frak u}} 
\nc{\gr}{\ensuremath{\frak r}} 
\nc{\gt}{\ensuremath{\frak t}}
\nc{\fff}{\ensuremath{\frak f}}
\nc{\gm}{\ensuremath{\mathfrak{mcf}}}
\nc{\gge}{\ensuremath{\mathfrak e}}
\nc{\cfupro}{\ensuremath{\cf(\upro)}}
\nc{\cfvpro}{\ensuremath{\cf(\vpro)}}
\nc{\gp}{\ensuremath{\frak p}}
\nc{\gk}{\ensuremath{\frak k}}

\nc{\add}[1]{\mbox{\ensuremath{{\rm add}(#1)}}}
\nc{\cov}[1]{\mbox{\ensuremath{{\rm cov}(#1)}}}
\nc{\unif}[1]{\mbox{\ensuremath{{\rm unif}(#1)}}}
\nc{\cof}[1]{{\mbox{\ensuremath{\rm cof}(#1)}}}

\nc{\addd}[2]{\mbox{\ensuremath{{\rm add}^{#1}(#2)}}}   
\nc{\covv}[2]{\mbox{\ensuremath{{\rm cov}^{#1}(#2)}}}   
\nc{\uniff}[2]{\mbox{\ensuremath{{\rm unif}^{#1}(#2)}}} 
\nc{\coff}[2]{{\mbox{\ensuremath{\rm cof}^{#1}(#2)}}}

\nc{\cd}{Cicho\'n's Diagram}

\nc{\MA}{\mbox{\rm MA}}
\nc{\PFA}{\mbox{\rm PFA}}
\nc{\OCA}{\mbox{\rm OCA}}
\nc{\GCH}{\mbox{\rm GCH}}
\nc{\CH}{\mbox{\rm CH}}
\nc{\zfc}{\mbox{\rm ZFC}}
\nc{\sch}{\mbox{\rm SCH}} 
\nc{\ZF}{\mbox{\rm ZF}}
\nc{\NCF}{\mbox{\rm NCF}} 
\nc{\FD}{\mbox{\rm FD}}   

\nc{\fourG}{\mbox{\rm 4G}} 
\nc{\fourI}{\mbox{\rm 4I}}   

\nc{\Borelhood}{Borel measurability} 
\nc{\Pieinseins}{\mbox{${\bf \Pi}^1_1$}}
\nc{\seinseins}{\mbox{${\bf\Sigma}^1_1$}}
\nc{\seinszwei}{\mbox{${\bf\Sigma}^1_2$}}
\nc{\seinsdrei}{\mbox{${\bf\Sigma}^1_3$}}
\nc{\Deleinszwei}{\mbox{${\bf\Delta}^1_2$}}

\nc{\up}{\ensuremath{{\cal U}\mbox{\ensuremath{\rm -prod}}\,\omega}}
\nc{\upp}{\ensuremath{{\cal U}'\mbox{\ensuremath{\rm -prod}}\,\omega}}
\nc{\upro}{\ensuremath{{\cal U}\mbox{\ensuremath{\rm -prod}}\,\om}}
\nc{\fupro}{\ensuremath{f({\cal U})\mbox{\ensuremath{\rm -prod}}\,\om}}
\nc{\vpro}{\ensuremath{{\cal V}\mbox{\ensuremath{\rm -prod}}\,\om}}
\nc{\fpro}{\ensuremath{{\cal F}\mbox{\ensuremath{\rm -prod}}\,\om}}

\nc{\cff}[1]{{\text{cf}\,(#1)}}           
\nc{\cu}{\ensuremath{\cal U}}             
\nc{\ai}{\ensuremath{\forall^\infty}}     
\nc{\ei}{\ensuremath{\exists^\infty}}     
\nc{\ww}{\ensuremath{\omega^\omega}}      

\nc{\gw}{groupwise dense}

\nc{\kk}{car\-dinal cha\-rac\-teris\-tic}
\nc{\joker}{\ast}
\nc{\gtc}{Galois-Tukey connection} 

\nc{\av}[1]{{\rm Av}_{#1}}
\nc{\eps}{\varepsilon}
\nc{\n}{{\bf n}}                 
\nc{\m}{{\bf m}}

\nc{\marginparr}[1]{}
\nc{\footnoteee}{} 
\nc{\footnotee}{}  

\newcommand{\cal}{\mathcal}

\nc{\divs}{{c_0 \setminus \ell^1}}
\nc{\divser}{(\divs, \leq^*)/\thickapproy}
\nc{\bfin}{\RO(\pfin \setminus\{0\},\subseteq^*)}
\nc{\bdivser}{\RO(\divser)}
\nc{\inc}{{\rm INC}}
\nc{\com}{{\rm COM}}
\nc{\thickapproy}{\makebox{}\!\!\thickapprox}
\nc{\approy}{\makebox{}\!\!\approx}
\nc{\lessi}{\leqslant}
\nc{\gessi}{\geqslant}
\nc{\interior}[1]{{\rm int}(#1)}
\nc{\closure}[1]{{\rm cl}(#1)}
\nc{\Vo}{Vojt\'a\v{s}}

\nc{\precedeseq}{\leq^*} 
\nc{\precedes}{\prec}
\nc{\stronger}{\leqslant_{\bf P}}
\nc{\underlline}[1]{\hat{#1}}
\nc{\PO}{{\bf P}}
\nc{\charak}{\text{ch}}
\nc{\symom}{{\rm{Sym}(\omega)}}
\begin{document}


\title{The relative consistency of $\gro < \cf(\symom)$}

\author{Heike Mildenberger and Saharon Shelah}
\thanks{The first author was supported by
a Minerva fellowship.}

\thanks{The second author's research 
was partially supported by the ``Israel Science
Foundation'', founded by the Israel Academy of Science and Humanities.
This is  the second author's work number 731}

\address{Heike Mildenberger,
Saharon Shelah,
Institute of Mathematics,
The Hebrew University of Jerusalem,
 Givat Ram,
91904 Jerusalem, Israel
}

\email{heike@math.huji.ac.il}
 
\email{shelah@math.huji.ac.il}

\begin{abstract}
We prove the consistency result from the title.
By forcing we construct a model of
$\gro= \aleph_1$, $\gb =  \cf(\symom) = \aleph_2$.
\end{abstract}

\subjclass{03E15, 03E17, 03E35}

\date{24.4.2000}

\maketitle

\setcounter{section}{-1}
\section{Introduction}\label{S0}

We recall the definitions of the three cardinal characteristics in
the title and the abstract. We write $A \subseteq^\ast B$ if
$A \setminus B$ is finite. We write $f \leq^\ast g$ if
 $f,g \in {}^\omega \omega $ and $\{ n \such f(n) > g(n) \}$ is
finite.

\begin{definition}\label{0.1}
\begin{myrules}
\item[(1)]
A subset $\cal G$ of $\potinf$ is called  groupwise dense if

 -- for all $B \in {\cal G}$, $A \subseteq^\ast B$ we have that
 $A \in {\mathcal G}$ and

 -- for every partition $\{
[\pi_i,\pi_{i+1}) \such i
\in \om \}$ of $\omega$ into finite intervals
there is an infinite set $A$ such that
$\bigcup\{[\pi_i,\pi_{i+1}) \such i \in A \} \in {\cal G}$.

The groupwise density number, \gro, is the smallest number of 
groupwise dense families with empty intersection.

\item[(2)] $\symom$ is the group of all
permutations of $\omega$.
If
$\symom = \bigcup_{i< \kappa} K_i$ and
$\kappa = \cf(\kappa)>\aleph_0$, $\langle K_i \such i < \kappa \rangle$ 
is increasing and continuous, $K_i$ is
a proper subgroup of $\symom$, we call
 $\langle K_i \such i < \kappa \rangle$ 
a cofinality witness. We call the minimal such $\kappa$ 
the cofinality of the symmetric group, short $\cf(\symom)$.

\item[(3)] 
The bounding number $\gb$ is
$$\gb = \min\{ |{\mathcal F}| \such {\mathcal F} \subseteq
{}^\omega \omega \,\wedge\, (\forall 
g \in {}^\omega \omega)(\exists f \in {\mathcal F})
f \not\leq^\ast g \}.$$

\end{myrules}
\end{definition}

Simon Thomas asked  whether $\gro \neq \cf(\symom)$
is consistent \cite[Question 3.1]{thomas:gd}.
In this paper we prove:

\begin{theorem}\label{0.2}
$\gro < \cf(\rm{Sym}(\omega))$ is consistent relative to $\zfc$.
\end{theorem}

\section{Forcings destroying many cofinality witnesses}\label{S1}

In this section we introduce two families of forcings that
will be used in certain steps of our planned iteration of
length $\aleph_2$. The plot is: If $\gb$ is large, there is some
way to destroy all shorter cofinality witnesses because 
by Claims~\ref{1.6} and \ref{1.5} none of the 
subgroups in a cofinality witness
contains all permutations respecting a given equivalence relation.
In our intended construction, we shall
 extend suitable intermediate models
 with a forcing built upon such an equivalence relation and
thus prevent possible cofinality
 witnesses to be lifted to the forcing extension and
all further extensions (Claim~\ref{1.4}).

Here we show some details
about destroying one cofinality witness that can be put separately
before we launch into an iteration. 
The additional task, to increase the bounding number along the way,
will be taken care of only in the next section.

\begin{definition}\label{1.1}
\begin{myrules}
\item[(1)]
We work with the following set of equivalence relations:
\begin{equation*}
\begin{split}
\makebox{} {\mathcal E}_{con}=
\{ E \such & E \mbox{ is an equivalence relation of } \omega,\\
& \mbox{each equivalence class is a finite interval
and }\\ 
& \omega = \liminf \langle |n/E| \such n < \omega \rangle \}.
\end{split}
\end{equation*}
We say $b \subseteq \omega$ respects $E\in {\mathcal E}_{con}$
if $(nEm \wedge m \in b) \rightarrow n \in b$. A partial permutation $\pi$ of
$\omega$ respects $E$ if $\dom(\pi)$ respects $E$ and
we have that  $n \in \dom(\pi) \rightarrow nE\pi(n)$.

\item[(2)]
Let $Q$ be the set of $p$ such that

\begin{myrules}
\item[(a)] $p$ is a permutation of some subset $\dom(p)$ of $\omega$,
\item[(b)] $\omega \setminus \dom(p)$ is infinite.
\end{myrules}
We order $Q$ by inclusion.

\item[(3)]
For $E \in {\mathcal E}_{con}$, $Q_E$ is the set of $p$ satisfying
(2)(a) -- (b) and additionally
\begin{myrules}
\item[(c)]  $p$ respects $E$.
\end{myrules}
\end{myrules}
\end{definition}

Part (1) of the following claim is important for later use, whereas 
part (2) will never be used directly.

\begin{claim}\label{1.2}
\begin{myrules}
\item[(1)] If $E \in {\mathcal E}_{con}$ and 
$p \in Q_E$ and $\name{\tau}$ is a  $Q_E$-name of an 
ordinal and $b$ is a finite subset of 
$\omega\setminus \dom(p)$ respecting $E$, then there is some $q$ such that 

\begin{myrules}
\item[(a)]
 $p \leq q$ and $b \subseteq \omega \setminus \dom (q)$,
\item[(b)]
if $\pi$ is a permutation of 
$b$ and it respects $E$ then $q \cup \pi$ forces a value to $\name{\tau}$.
\end{myrules}

\item[(2)] $Q_E$ is proper, ${}^\omega \omega$-bounding, nep
(see \cite{Sh:630}) and Souslin.
\end{myrules}
\end{claim}
\proof (1) Note that there are only finitely many permutations of $b$
(that respect $E$).
So we can treat them consecutively and find stonger and stronger $q$'s.

(2) Let $N \prec H(\chi,\in)$ be such that $Q_E \in N $ and
$p \in N$, $\chi \geq (2^\omega)^+$. Let $\name{\tau_n}$, $n \in \omega$, be a list of all
$Q_E$-names for ordinals that are in $N$. Let
$b_n$, $n \in \omega$, be a list of pairwise disjoint
 $E$-classes such that
$\bigcup_{n \in \omega } b_n $ is infinite. 
Now take $q_n$ by induction starting with $q_0 = p$. If $q_n$ is chosen,
take $i(n)$ 
such that $\dom(q_n) \cap b_{i(n)} = \emptyset$.
Now take $q_{n+1}$ treating $q_n$, $\name{\tau_n}$ and $b_{i(n)} $ as in 
the proof of part (1). We have that $q = \bigcup q_n \in Q_E$ and that
$q \Vdash_{Q_E} (\forall n \in \omega) \name{\tau_n} \in \check{N}$.
By \cite[III, Theorem 2.12]{Sh:h}, $Q_E$ is proper.

$Q_E$ is ${}^\omega \omega$-bounding: 
Let $\name{f}$ be a name for a function from
$\omega$ to $\omega$.
Again let $b_n$, $n \in \omega$, be a list of pairwise disjoint
 $E$-classes such that
$\bigcup_{n \in \omega } b_n $ is infinite. 
Now take $q_n$ by induction starting with
 $q_0 = p$. If $q_n$ is chosen,
take $i(n)$ 
such that $\dom(q_n) \cap b_{i(n)} = \emptyset$.
Now take $q_{n+1}$ treating $q_n$, $\name{\tau_n}$ and $b_{i(n)} $ as in 
part (2) of this claim and look which values for $\name{f}(n)$ 
the finitely many permutations in (1)(b) force. Take $g(n)$ 
to be the maximum of them. We have that $q = \bigcup q_n \in Q_E$ and that
$q \Vdash_{Q_E} (\forall n)  \name{f}(n) \leq g(n)$.

nep (non-elementary properness): 
We use much less than $N \prec H(\chi,\in)$. We use that $E\in N
\subseteq H(\chi,\in)$.
See \cite{Sh:630}.

Souslin: $p \in Q_E$, $q \leq q$ and $p \perp q$ can be expressed in  
$\Sigma^1_1(E)$-formulas.
\proofend

We shall work
with the following special subsets of $\symom$.

\begin{definition}\label{1.3}
\begin{myrules}
\item[(1)]
For $E \in {\mathcal E}_{con}$ and $A \subseteq \omega$ we define:
$$ S_{E,A}:=
\{ \pi \in Q_E \such \pi \restriction (\omega \setminus A) = id \}.$$
\item[(2)]
We set ${\mathcal F} :=
\{ f \such f \in {}^\omega \omega, f(n) \geq n,
\lim\langle f(n) - n \such n \in \omega \rangle = \infty \}$.
For $f \in {\mathcal F}$ we set
$ S_{f} :=
\{ \pi \in \symom \such 
(\forall n) \pi(n) \leq f(n) \}.$
\end{myrules}
\end{definition}

The following claim describes the basic step in order to increase
$\cf(\symom)$.

\begin{claim}\label{1.4}
Assume
\begin{myrules}
\item[(a)] 
$\langle K_i \such i < \kappa \rangle$ 
is a cofinality witness,
\item[(b)] $\name{R}$ is a $Q_E$-name of a forcing notion,
\item[(c)] $E \in {\mathcal E}_{con}$, and for no $i < \kappa$ and
coinfinite  $A
\in [\omega]^\omega$ respecting $E$ we have that $K_i \supseteq S_{E,A}$.
\end{myrules}
Then in ${\bf V}^{Q_E \ast \name{R}}$ we cannot find a cofinality witness
$\langle K'_i \such i < \kappa \rangle$ such that $\bigwedge_{i< \kappa}
\left(K'_i \cap \symom^{\bf V} = K_i\right)$.
\end{claim}

\proof
Let $\name{f}= \bigcup\{ p \such p \in \name{G_{Q_E}} \}$ be a $Q_E$-name 
of a permutation of $\omega$.
It suffices that 
\begin{equation}\tag{$\ast$} \label{ast}
\begin{split}
\Vdash_{Q_E} &\mbox{``for unboundedly many $i < \kappa$,}\\
&\mbox{for some $g \in K_i$
we have $\name{f} \circ g 
\in K_{i+1} \setminus K_i$.''}
\end{split}
\end{equation}

Why does this suffice? Suppose that \eqref{ast} holds and
we had found a cofinality witness
$\langle K'_i \such i < \kappa \rangle$ 
in ${\bf V}^{Q_E \ast \name{R}}$ such that $\bigwedge_{i< \kappa}
\left(K'_i \cap \symom^{\bf V} = K_i\right)$.
Let $G$ be $Q_E \ast \name{R}$-generic over ${\bf V}$.
Take $j < \kappa$ such that $\name{f}[G] \in K'_j$. Then 
we find according to
\eqref{ast} some $i \geq j$ and some $g \in K_i$ such that
$\name{f[G]}\circ g \in K_{i+1} \setminus K_i \subseteq {\bf V}$.
But this contradicts the facts that $\name{f}[G] \circ g \in K'_{i}$
(because this is a subgroup) and
$K'_i \cap \symom^{\bf V} = K_i$.

\smallskip

Proof of \eqref{ast}: Let $p \in Q_E$ and $j < \kappa$.
Let $\omega \setminus \dom(p)$ be the disjoint union of $A_0,A_1$, 
both infinite subsets of $\omega$ respecting $E$.
Let $g_0\in \symom$ be such that 
$\{n \such g_0(n) \neq n \} = A_0$. Let $g_0 \in
 K_{i(\ast)}$, $i(\ast) > j$. 
By assumption $S_{E,A_0}$ is not included in any $K_i$, so in
particular not included in $K_{i(\ast)}$.
Hence there is $g_1 \in S_{E,A_0} \setminus K_{i(\ast)}$.
Take $i$ such that $g_1 \in K_{i+1} \setminus K_i$. 
Necessarily we have $\kappa > i \geq i(\ast) > j$. 
Now there is a permutation $f$ of $A_0$ respecting $E$ 
such that $f$ is an isomorphism from
$(A_0,g_1)$ onto $(A_0, g_0)$.
Namely set $f(g_0(n)) = g_1(n)$.
 Hence
$n \in A_0 \Rightarrow f(g_0(n)) = g_1(n)$. Let
$q = p \cup f$. The condition
$q$ forces that $\name{f} \circ g_0  = g_1$, $g_1 
\in K_{i+1} \setminus K_i$, and $i \in (j,\kappa)$, 
$g_0 \in K_{i(\ast)} \subseteq K_i$, so 
\eqref{ast} is proved.
\proofend

\begin{claim}\label{1.5}
Assume that $\langle K_i \such i < \kappa \rangle$ is a cofinality witness.
Assume that $K_0$ contains all permutations that move only finitely 
many points.
Then the following are equivalent:

\begin{myrules}
\item[$(\alpha)$]
For some $E \in {\mathcal E}_{con}$, 
for no $i < \kappa$, coinfinite $A \in [\omega]^{\aleph_0}$ 
we do have $K_i \supseteq S_{E,A}$.
\item[($\beta$)] 
For every $E \in {\mathcal E}_{con}$, 
for no $i < \kappa$, coinfinite $A \in [\omega]^{\aleph_0}$ 
we do have $K_i \supseteq S_{E,A}$.
\item[($\gamma$)] 
For some $f \in {\mathcal F}$, 
for no $i < \kappa$ do we have that $S_{f} \subseteq K_i$.
\item[($\delta$)] 
For every $f \in {\mathcal F}$,
for no $i < \kappa$ do we have that $S_{f} \subseteq K_i$.
\end{myrules}
\end{claim}

\proof
The implications $(\beta) \Rightarrow (\alpha)$
 and  $(\delta) \Rightarrow (\gamma)$ are 
 trivial. We shall not use $(\beta) \Rightarrow (\alpha)$ but
close a circle
of implications as follows:
$(\beta) \Rightarrow (\delta)$ and $(\alpha) \Rightarrow (\beta)$
and $(\gamma) \Rightarrow (\alpha)$.

Now we prove
$\neg (\delta) \Rightarrow \neg(\beta)$.
Let $f$ and $i^\ast$ exemplify the failure of $(\delta)$.
\nothing{Let $E \in {\mathcal E}_{con}$.
Let $\langle n_k \such k \in \omega \rangle$ be such that $n_0 = 0$,
$n_k < n_{k+1}$ and each interval $[n_k,n_{k+1})$ is an $E$-equivalence 
class. We can choose by induction on $k \in \omega$ some 
$m_k$ such that 
\begin{myrules}
\item[(i)] $(\forall j < k) (m_j + n_{j+1} < m_k < \omega)$.
\item[(ii)] $(\forall m) (m_k \leq m \rightarrow f(m) \geq m + n_{k+1})$.
(So, here we use the definition of $S_{f}$.)
\end{myrules}
Let $\ell \in \{0,1\}$ and $A_\ell = \bigcup \{ [ n_{2k +\ell},
n_{2k + \ell +1}) \such k < \omega \}$, so $A_0, A_1$ is a partition 
of $\omega$ to infinite $E$-respecting sets. Now for $\ell = 0,1$ let
$g_\ell$ be a permutation of $\omega$ which maps for
every $k \in \omega$ the interval 
$[n_{2k + \ell},n_{2k + \ell + 1 })$ into the interval 
$[m_{2k + \ell},m_{2k + \ell} + n_{2k + \ell +1} - n_{2k + \ell})$. 
Hence multiplication  by
$g_\ell$ maps $S_{E,A_\ell}$ into
\begin{equation*}
\begin{split}
B_\ell =&\{ \pi \in \symom \such \pi \mbox{ maps } 
[n_{2k+\ell}, n_{2k+\ell+1})\\
& \mbox{  onto }
[m_{2k+\ell}, m_{2k+\ell} + n_{2k+ \ell +1} - n_{2k + \ell})\\ 
& \mbox{ for $k < \omega$ and is the identity otherwise$\},$}
\end{split}
\end{equation*}
which is included in $\langle S_{f} \cup K_0\rangle$ (where that
angles denote the subgroup of $\symom$ that is generated by
the set in  between then) 
by the choice of the $m_i$'s.
Indeed, for $\ell = 0,1$ for $r \in [m_{2k +\ell},
m_{2k+ \ell + 1 })$, $\pi \in B_\ell$ we have that
$\pi(r) \leq 
m_{2k + \ell} + n_{2k + \ell + 1}
-n_{2k + \ell} \leq f(r)$ for almost all $r\in \omega$.
The finitely many exceptions are taken care of by $K_0$.

 Hence any subgroup
of $\symom$ which includes $S_{E,f} \cup \{g_0,g_1\}
\cup K_0$ includes 
$S_{E,A_0} \cup S_{E,A_1}$. For some $j(\ast) \in
[i(\ast),\kappa)$ we have that $g_0,g_1 \in K_{j(\ast)}$, and hence
by $\neg(\delta)$,
$S_{E,A_0} \cup S_{E,A_1} \subseteq K_{j(\ast)}$. 
So for every $E$ there are $j(\ast)<\kappa$ and some infinite, coinfinite
$A$ that exemplify the failure of $(\alpha)$.}

By the definition of ${\mathcal F}$ we have that
$\lim\langle f(n) - n \such n \in \omega \rangle = \infty$.
Hence we may choose a strictily increasing sequence $
\langle k_i \such i \in \omega \rangle$ such that $(\forall
i \in \omega) (\forall n \geq k_i) (f(n) \geq i + n)$.
Then we take $E = \{ [k_i, k_i + i) \such i \in \omega\}
\cup \{[k_i +i, k_{i+1}) \such i \in \omega\}$ and
$A = \bigcup_{i \in \omega} [k_i,k_{i+1})$. $A$ is infinite and
coinfinite.
Then we have that $S_{E,A} \subseteq S_f \subseteq K_{i^\ast}$,
so $\neg(\beta)$.

\smallskip

Now we show $\neg(\beta)$ implies $\neg(\alpha)$.
This follows from 

\begin{claim*}
For all $E,E' \in {\mathcal E}_{con}$ 
 there are $f_1,f_2 \in \symom$ such that
for any $A \subseteq \omega$ we have
$$S_{E,A} \subseteq (f_1 \circ S_{E', f_1^{-1}[A]} \circ f_1^{-1})\circ
(f_2 \circ S_{E', f_2^{-1}[A]}\circ f_2^{-1}).$$
\end{claim*}

\proof Enumerate the $E$-classes with order type $\omega$.
Let $f_1$ inject the even-numbered $E$-classes into high enough
(there are large enough ones by the definition of ${\mathcal E}_{con}$)
$E'$ classes. The $E'$-classes need not be covered, it is enough that
$nEm \rightarrow f_1(n) E' f_1(m)$. We fill this function up to
a permutation of $\omega$ and call it $f_1$.
Let $f_2$ do the same with the odd-numbered $E$-classes.
If $g \in S_{E,A} $ then $g = g_1 \circ g_2$ where $g_1$ is the identity on
odd-numbered $E$-classes and $g_2$ is the identity on even-numbered
$E$-classes. We have that $f_i^{-1}\circ
 g_i\circ f_i \in S_{E,f_i^{-1}[A]}$ for $i=1,2$
and thus the claim is proved.

\smallskip

To complete a cycle of implications, we show
$\neg (\alpha) \Rightarrow \neg(\gamma)$.
\nothing{Let $\langle n_k \such k \in \omega \rangle$ be such that $n_0 = 0$,
$n_k < n_{k+1}$, $n_k$ even 
 and each interval $[n_k,n_{k+1})$ is an $E$-equivalence 
class.}
To prove $\neg(\gamma)$ let $f \in {\mathcal F}$. 
We choose by induction on $k \in \omega$, $m_i$ such that $m_0=0$,
$m_{k+1} > m_k$ and 
$(\forall n< m_k) (f(n) <m_{k+1})$.
Now we take $n_i$ by induction on $i$ such that 
$n_0=0$, $n_{i+1} > n_i$ and 
$(\forall m \geq m_{n_{i+1}}) (\pi(m) \geq m_{n_i})$.

Now we define two equivalence relations.

\begin{equation*}\begin{split}
E_0 &= \{ [m_{n_k},m_{n_{k+2}}) \such k \in \omega \},\\
E_1 &= \{ [m_{n_{k+1}}, m_{n_{k+3}}) \such
k \in \omega \} \cup \{ [0, m_1)\}.
\end{split}
\end{equation*}

For $\mu \in \{0,1,2,3\}$ let $A_\mu = \bigcup \{ [ m_{n_k},
m_{n_{k+3}}) \such k < \omega, k = \mu \text{mod} 4 \}$.

Now note that 
\begin{myrules}
\item[$(\ast)_1$]

If $\pi \in S_{f}$ then we can find $\pi_\ell \in
S_{E_\ell,\omega}$ for $\ell = 0,1$ such that $\pi = 
\pi_1 \circ \pi_0$. Why?

By our choice of $\langle m_i \such i \in \omega \rangle$
and $\langle n_i \such i \in \omega \rangle$ and 
 $E_\ell$, for any $x \in \omega$, 
$x E_0 \pi(x)$ or $x E_1 \pi(x)$. Now we choose
$\pi_0(x)$ and $\pi_1(x)$ by cases.

If $x $ and $\pi(x)$ are in the same $E_0$-class
and $\pi(x) E_0 \pi(\pi(x))$, then
we set $\pi_0(x) = \pi(x)$ and $\pi_1(\pi(x))= \pi(x)$.
So we have $\pi(x) = \pi_1 \circ \pi_0(x)$.

If $x $ and $\pi(x)$ are in the same $E_0$-class
and not $\pi(x) E_0 \pi(\pi(x))$, then
we set $\pi_0(x) = y$ and $\pi_1(y) = \pi(x)$ for some
$y E_0 x$ such that $\pi(y) \neq y$ and $y E_1 \pi(x)$
and $y E_0 \pi(y)$.
(If there are not enough such $y$, just take the classes
of ``double width''. We also assume w.l.o.g.
that $\pi$ has no fixed points.) 
Then we have that $\pi_0$ respects $E_0$ in the point $x$, and
$\pi_1$ respects $E_1$ in the point $y$ and $\pi(x) =\pi_1 \circ
\pi_0(x)$.

If $x$ and $\pi(x)$ are not in
the same $E_0$-class, then we have
that $x E_1 \pi(x)$. If not
$\pi^{-1}(x) E_0 x$ then we set $\pi_1(x)=\pi(x)$ and $\pi_0(x) = x$.

If $x$ and $\pi(x)$ are not in
the same $E_0$-class, then we have
that $x E_1 \pi(x)$. If 
$\pi^{-1}(x) E_0 x$ then we also set $\pi_0(x)=x$
and $\pi_1(x) = \pi(x)$.
Note that the pair $(\pi^{-1}(x),x)$ falls under the second case and that
hence there is no conflict in our settings, i.e.\
also $\pi_0$ and $\pi_1$ can be chosen as permutations.

Then we have that $\pi = \pi_1 \circ \pi_0$.

\item[$(\ast)_2$]
Let $\ell= 0,1$ and $\pi_\ell$ be as above.
Then we can find $\psi_{\ell,\mu} \in
S_{E_\ell,A_\mu}$ for $\mu = 0,1,2,3$ such that $\pi_\ell = 
\psi_{\ell,3} \circ \psi_{\ell,2} \circ \psi_{\ell,1} 
\circ \psi_{\ell,0}$. Why?
For all $x \in \omega$ there are three $\mu$'s such that
$x \in A_\mu$ and three $\mu$'s such that
$\pi_\ell(x) \in A_\mu$. Hence we can find $\mu$ (indeed, 
two $\mu$'s) such that
$x,\pi_\ell(x) \in A_\mu$, and such we may chose some
$\psi_{\ell,\mu} \in S_{E,A_\mu}$ such that $\pi_\ell(x) = 
\psi_{\ell,\mu}(x)$ and
such that $\psi_{\ell,\mu}$ restricted to $\omega \setminus A_\mu$
 is the identity and such that 
$\pi_\ell = 
\psi_{\ell,3} \circ \psi_{\ell,2} \circ \psi_{\ell,1} 
\circ \psi_{\ell,0}$. 

\item[$(\ast)_3$] Let for $\ell = 0,1$ the infinite, coinfinite
set $A^\ell$ and the ordinal $i^\ell(\ast)$ be
as in $\neg(\alpha)$ for $E_\ell$. For $\mu < 4$ there is 
$g_{\mu} \in \symom$ mapping
$\omega \setminus A_\mu$ into $\omega \setminus A^\ell$
such that $(\forall k_0,k_1 \in \omega \setminus A_\mu)
(k_0 E_\ell k_1 \Leftrightarrow g_\mu(k_0) E_\ell g_\mu(k_1))$,
hence for $\ell = 0,1$
conjugation by $g_\mu$ maps $S_{E_\ell,A_\mu}$ into $S_{E_\ell,A^\ell}
\subseteq K_i$. 
\end{myrules}

By our assumption $\neg(\alpha)$ 
we have some $i^\ell(\ast) \in \kappa$ such that
$S_{E_\ell,A^\ell} \subseteq K_{i^\ell(\ast)}$ for $\ell = 0,1$
and $\mu=0,1,2,3$. Let $i(\ast)= \max(i^0(\ast),i^1(\ast))$.
For some $j(\ast) \in
[i(\ast),\kappa)$ we have that $g_\mu \in K_{j(\ast)}$
for $\mu=0,1,2,3$, and 
$S_{E_\ell,A_\mu} = g_\mu \circ S_{E_\ell,A^\ell} \circ g^{-1}_\mu 
\subseteq K_{j(\ast)}$, hence
$S_{f} \subseteq K_{j(\ast)}$, that is, $\neg(\gamma)$.
\proofend

\begin{claim}\label{1.6}
Assume that
$\langle K_i \such i < \kappa \rangle$ is a cofinality witness
such that $K_0$ contains all the permutations that
move only finitely any points.
If $\gb >\kappa$, then clause $(\gamma)$  of Claim~\ref{1.5} holds 
(and hence all the other clauses hold as well).
\end{claim}
\proof For each $ i < \kappa$ choose $\pi_i 
\in \symom \setminus K_i$. Since $\gb > \kappa$ 
there is some $f \in {}^\omega \omega$ such that 
$(\forall i < \kappa) (\forall^\infty n) (\pi_i(n) < f(n))$ and w.l.o.g.\
$f \in {\mathcal F}$. 
if $S_f$ were a subset of $K_i$, then we had that $\pi_i \in K_i$, which is not the case.
So $f$ exemplifies clause $(\gamma)$ of Claim~\ref{1.5}.
\proofend

\begin{definition}\label{1.7}
\begin{myrules}
\item[(1)] Let $E\in {\mathcal E}_{con}$.
We set 
\begin{equation*}
\begin{split}
Q'_E = \{ f \such & \mbox{$f$ is a permutation of some coinfinite subset
of $\omega$ such that }\\
&\mbox{(a) } n \in \dom(f) \Rightarrow n E f(n),\\
&\mbox{(b) for every $k < \omega$ for some $n$ we have }
k \leq |(n/E) \setminus \dom(f)| \}.
\end{split}
\end{equation*}
The order is by inclusion.

\item[(2)] We call $\bar{f} = \langle f_i \such i < \alpha \rangle$,
$Q'_E$-o.k.\ if $\alpha \leq \omega_1$ and for $i\leq j<\alpha$, 
$f_i \subseteq^\ast f_j\in Q'_E$ (i.e. 
$\{ n \in \dom(f_i) \such n \not\in \dom(f_j) \vee 
f_i(n) \neq f_j(n) \}$ is finite).
For $\bar{f}$ being $Q'_E$-o.k.\ we set
$Q'_E(\bar{f})= \{ g \in Q_E' \such g =^\ast f_i $ for some $i\}$,
where $f_i =^\ast g$ iff $f_i \subseteq^\ast g$ and $g \subseteq^\ast f_i$.
The order is inherited from $Q'_E$.
\end{myrules}

\end{definition}
 
\begin{remarks*}
1)
Claims~\ref{1.4} and \ref{1.5} hold for $Q'_E$ as well
with the analogously modified definition of
$S'_{E,A}$.  This is 
shown with the same proofs. The domains of the involved partial permutations
 must be arranged 
such that they respect \ref{1.7}(1)(b), but they need not be unions of
equivalence classes.
The $q \in Q_E $ fulfil requirement \ref{1.7}(1)(b) automatically, 
because we have that 
 $\lim\langle |n/E|
\such n \in \omega \rangle = \omega$ and that the domain
of $q$  needs to be 
coinfinite and needs to be a union of equivalence classes.

2) Both $Q_E$ and $Q'_E$ can serve for our purpose. 
$Q'_E$ exhibits the following ``independence of $E$'':
For $E_0,E_1 \in {\mathcal E}_{con}$ $(\forall p \in Q'_{E_1})$
$(\exists q)$ $(p \leq q \in Q'_{E_1} \wedge (Q'_{E_1})_{\geq p}
\cong Q'_{E_0})$.

3)
Note that for $\alpha< \omega_1$, if  $\bar{f} =
\langle f_\beta \such \beta \in \alpha 
\rangle$ $Q'_E$-o.k., then we have that $Q'_E(\bar{f})$  is Cohen forcing.
\end{remarks*}
 
\begin{claim}
\label{1.8}
Let $E$ be as in Definition~\ref{1.7}. 
\begin{myrules}
\item[(1)] $Q'_E$ is proper, even strongly proper, with the Sacks property
(the last is more than $Q_E$).
\item[(2)]
If $p \in Q'_E$ and a sequence 
$\langle w_n \such n \in \omega \rangle$ of pairwise disjoint
finite subsets of $\omega$ are given, then
we find an infinite $u \subseteq \omega$ such that
$\langle w_n \such n \in u \rangle$
and 
$(\forall n) (\exists m) (w_n \subseteq m/E) $ and $w_n \cap \dom(p)
= \emptyset$ and $n_1 < n_2 \Rightarrow \forall m_1 \in w_{n_1}
\forall m_2 \in w_{n_2} \neg m_1 E m_2$, and for
 every permutation $f$ of
$\bigcup w_n$ which respects $E$ we have that $p \leq p \cup f 
\in Q'_E$.

\item[(3)]
If $\bar{f} $ is as in \ref{1.7}(2), and $\alpha < \omega_1$ and 
$Q'_E(\bar{f})\subseteq M$, $\omega +1 \subseteq M$ $ \subseteq 
\mbox{$(H(\chi),\in)$}$, 
$M$ a countable model of $\zfc^-$, then we can find $f_\alpha$ such that 
\begin{myrules}
\item[(a)] $\bar{f} \concat f_\alpha$ is $Q'_E$-o.k.
\item[(b)] 
If $\bar{f} \concat f_\alpha \triangleleft \bar{f'} $ and $\bar{f'}$ is
$Q'_E$-o.k., then $f_\alpha$ is $(M,Q'_E(\bar{f'}))$-generic.
In fact, for every predense $I 
\subseteq Q'_E(\bar{f'})$ from $M$ some finite 
$J \subseteq I$
is predense above $f_\alpha$ in 
$Q'_E(\bar{f'})$. In fact, $J$ does not depend on $\bar{f'}$.
\end{myrules}
\end{myrules}
\end{claim}

\proof
(1) We prove the Sacks property.
Let $\name{f} \in V^{Q'_E} \cap {}^\omega \omega$.
 We take $b_n$ as in the proof of the ${}^\omega \omega$-boundedness
 for $Q_E$ (which applies also to $Q'_E$)
in Claim~\ref{1.2}, but we do not require that
$b_n$ respects $E$. Additionally we choose 
$b_n$ so small that there are only fewer than $n$ permutations of $b_n$.
Then we take $q_n$ as there and collect into $S(n)$ all the possible
values forced by $q_n \cup \pi$ for $\name{f}(n)$, when $\pi$ ranges 
over the permutations of $b_n$.

(2) Easy.

(3) Let $\langle \bar{f'}^n \such n \in \omega \rangle$ enumerate all the
$\alpha \leq \omega_1$-sequences in $M$ that are $Q'_E$-o.k.
Let $\name{\tau_n}$, $b_n$, $n \in \omega$ be as in the proof of \ref{1.2}.
Now we choose $f_\alpha^n$ $\subseteq^\ast$-increasing with $n$,
and  $i(n)$ strictly increasing with $n$ such that
$b_{i(n)} \cap \dom(f_\alpha^n)= \emptyset$ and such that 
if $\bar{f}\concat f^n_\alpha \triangleleft \bar{f'}^n $ then
$f_\alpha^n \Vdash_{Q'_E} \name{\tau_n} \in V$. This is done
with the finitely many permutations of a suitable $b_{i(n)}$ as in \ref{1.2}.
Note that 
$f_\alpha^n \Vdash_{Q'_E} \name{\tau_n} \in V$
and $\bar{f} \concat f_\alpha^n \triangleleft \bar{f'}$ implies
$f_\alpha^n \Vdash_{Q'_E(\bar{f'})} \name{\tau_n} \in {\bf V}$,
independent of the choice of $\bar{f'}$. 
We set $f_\alpha = \bigcup_{n \in \omega} f^n_\alpha$,
and by one of the equivalent characterizations of
$(M,Q'_E(\bar{f'}))$-genericity \cite[III, Theorem 2.12]{Sh:h}
we are done.
\proofend

\nothing{
Finally we show the analogue of \ref{1.4}
also for the ``thin'' forcings $Q'_E(\bar{f'}))$
if the length of $\bar{f}$ is countable.
\begin{claim}\label{1.9}
Assume
\begin{myrules}
\item[(a)] 
$\langle K_i \such i < \kappa \rangle$ 
is a cofinality witness,
\item[(b)] $\name{R}$ is a $Q'_E(\bar{f})$-name of a forcing notion,
\item[(c)] $E \in {\mathcal E}_{con}$, and for no $i < \kappa$ and $A
\in [\omega]^\omega$ respecting $E$ we have that $K_i \supseteq S_{E,A}$.
\end{myrules}
Then in ${\bf V}^{Q'_E(\bar{f}) \ast \name{R}}$ we cannot find a cofinality witness
$\langle K'_i \such i < \kappa \rangle$ such that $\bigwedge_{i< \kappa}
\left(K'_i \cap \symom^V = K_i\right)$.
\end{claim}
\proof
Again, let $\name{f}= 
\bigcup\{ p \such p \in \name{G_{Q'_E(\bar{f})}} \}$ be a 
$Q'_E(\bar{f})$-name 
of a permutation of $\omega$.
It suffices that 
\begin{equation}\tag{$\ast$} \label{ast}
\begin{split}
\Vdash_{Q'_E(\bar{f})} &\mbox{``for unboundedly many $i < \kappa$,}\\
&\mbox{for some $g \in K_i$
we have $\name{f} \circ g 
\in K_{i+1} \setminus K_i$.''}
\end{split}
\end{equation}
\smallskip
Proof of \eqref{ast}: Let $p \in Q'_E(\bar{f})$ and $j < \kappa$.
Let $\omega \setminus \dom(p) $ be the disjoint union of $A_0,A_1$, 
both infinite subsets of $\omega$ respecting $E$ such that
additionally for all $\beta < \lgg(\bar{f})$
we have that $\dom(f_\beta) \cap A_i$ is finite. (This works only for
$\lgg(\bar{f})$ countable.)
Let $g_0\in \symom$ be such that 
$\{n \such g_0(n) \neq n \} = A_0$. Let $g_0 \in
 K_{i(\ast)}$, $i(\ast) > j$. 
By assumption $S_{E,A_0}$ is not included in any $K_i$, so in
particular not included in $K_{i(\ast)}$.
Hence there is $g_1 \in S_{E,A_0} \setminus K_{i(\ast)}$.
Take $i$ such that $g_1 \in K_{i+1} \setminus K_i$. 
Necessarily we have $\kappa > i \geq i(\ast) > j$. 
Now there is a permutation $f$ of $A_0$ respecting $E$ 
such that $f$ is an isomorphism from
$(A_0,g_1)$ onto $(A_0, g_0)$.
Namely set $f(g_0(n)) = g_1(n)$.
 Hence
$n \in A_0 \Rightarrow f_0(g_0(n)) = g_1(n)$. Let
$q = p \cup f$. The condition
$q$ is in $Q'_E(\bar{f})$ and forces that $\name{f} \circ g_0  = g_1$, $g_1 
\in K_{i+1} \setminus K_i$, and $i \in (j,\kappa)$, 
$g_0 \in K_{i(\ast)} \subseteq K_i$, so 
\eqref{ast} is proved.
\proofend
}

\section{Arranging $\gro= \aleph_1, \gb=\cf(\symom)=
\aleph_2$}
\label{S2}

Starting from a ground model with a suitable diamond sequence
we find a forcing extension with the constellation from the section headline.
The requirements on the ground model
 can be established by a well-known forcing 
(see \cite[Chapter 7]{Kunen})
starting from any ground model, and are also true in $L$   (see
\cite{Devlin}).

 \begin{definition}\label{2.1}
\begin{myrules}
\item[(1)] We say ${\mathcal A}$ is a $(\kappa,\gro)$-witness
if $\kappa = \cf(\kappa) > \aleph_0$ and 
\begin{myrules}
\item[$(\alpha)$] ${\mathcal A} \subseteq [\omega]^{\aleph_0}$,
\item[$(\beta)$] if $k < \omega$ and $f_\ell \colon \omega \to \omega$ 
is injective for $\ell < k$ then for some ${\mathcal A}' 
\subseteq {\mathcal A}$ of cardinality $< \kappa$ we have that for any 
$A$ that is a finite union of members of ${\mathcal A} \setminus 
{\mathcal A}'$ 
$$\{ n \such \bigwedge_{\ell < k} f_\ell(n) \not\in A \}
\mbox{ is infinite.}
$$
\end{myrules}

\item[(2)] We say $\bar{M}$ $\kappa$-exemplifies ${\mathcal A}$ if 
\begin{myrules}
\item[(a)] ${\mathcal A} $ is a $(\kappa,\gro)$-witness,
\item[(b)] $\bar{M}= \langle M_i \such i < \kappa \rangle$ is 
$\prec$-increasing and continuous, and $\omega + 1 \subseteq M_0$ and 
${\mathcal P}(\omega)
\subseteq \bigcup_{i < \kappa} M_i$,
\item[(c)] $M_i \subseteq (H(\chi),\in)$ is a model of 
$\zfc^-$ and $|M_i| < \kappa$ and $(M_i \models |X| < \kappa) 
\Rightarrow X \subseteq M_i$,
\item[(d)] $\bar{M} \restriction (i+1) \in M_{i+1}$,
\item[(e)] for $i$ non-limit, 
there is ${\mathcal A}_i \in M_i$ such that 
${\mathcal A} \cap M_i = {\mathcal A}_i$,
\item[(f)] if $i < \kappa$, $k < \omega$ and 
$f_\ell \in M_i$ is an injective function from $\omega$ to $\omega$ 
for $\ell < k$, and $k' < \omega$, $A_\ell \in 
{\mathcal A} \setminus M_i$ for $\ell < k'$, then 
$$\{ n \such \bigwedge_{\ell < k} f_\ell(n) \not\in
A_0 \cup \cdots \cup A_{k'-1} \} \mbox{ is infinite.}
$$
\end{myrules}
\nothing{\item[(3)] $M$ is called ${\mathcal A}$-nice if it 
satisfies the requirements on one of the $M_i$ above, e.g.\ on $M_0$.
}
\item[(3)] We say $\bar{M}$ leisurely exemplifies ${\mathcal A}$ if
(a) to (f) above are fulfilled and additionally;
\begin{myrules}
\item[(g)] $\kappa = \sup\{ i \such M_{i+1} \models
\mbox{``}{\mathcal A}_{i+1} = \aleph_0\mbox{''}\}.$
\end{myrules}
\nothing{
\item[(5)] $M$ is ${\mathcal A}$-good if it is ${\mathcal A}$-nice and 
$M\models |{\mathcal A}_M| = \aleph_0$ where 
${\mathcal A}_M \in M$ and ${\mathcal A} \cap M =
{\mathcal A}_M$.
}
\end{myrules}
\end{definition}

\nc{\pAA}{{\mathcal A}}
\nc{\nAA}{{\name{\mathcal A}}}

\begin{definition}\label{2.2}
\begin{myrules}
\item[(1)] We say $(P,\nAA)$ is a $(\mu,\kappa)$-approximation if

\begin{myrules}
\item[$(\alpha)$] $P$ is a c.c.c.\ forcing notion, $|P| \leq \mu$,
\item[$(\beta)$] $\nAA$ is a set of $P$-names of members 
of $([\omega]^{\aleph_0})^{{\bf V}^P}$, each hereditarily countable, 
and for simplicity
they are forced to be pairwise distinct,
\item[($\gamma$)] $\Vdash_P \mbox{``}\nAA 
\mbox{ is a $(\kappa,\gro)$-witness.''}$
\end{myrules}

\item[(2)] If $\mu= \kappa$ we may write just 
$\kappa$-approximation.
If $\kappa = \aleph_1$ we may omit it. 
We write $(\ast,\kappa)$-approximation if it is 
a $(\mu,\kappa)$-approximation for some $\mu$.
\item[(3)]
$(P,\nAA_1) \leq^\kappa_{app} (P_2,\nAA_2)$ 
if:
\begin{myrules}
\item[(a)]
$(P_\ell, \nAA_\ell)$ is a $(\ast,\kappa)$-approximation.
\item[(b)] $P_1 \lessdot P_2$,
\item[(c)] $\nAA_1 \subseteq \nAA_2$ (as a set of names, for simplicity),
\item[(d)] 
if $k < \omega$ and $\name{A_0}, \dots , \name{A_{k-1}} \in \nAA_2 \setminus
\nAA_1$ then 
\begin{equation*}
\begin{split}
\Vdash_{P_2} &\mbox{`` if $B \in ([\omega]^{\aleph_0})^{V^{P_1}}$},\\
& f_\ell \in ({}^B \omega)^{V^{P_1}} \mbox {for } \ell < k
\mbox{ are injective, then }\\
& \left\{ n \in B \such \bigwedge_{\ell <k} f_\ell(n) 
\not\in \bigcup_{\ell <k} 
\name{A_\ell} \right\} \mbox{ is infinite''.}
\end{split}
\end{equation*}
\end{myrules}
\end{myrules}
\end{definition}

\begin{remark*}
We mean
$\nAA_1 \subseteq \nAA_2$ as a set of names.
It is no real difference if $\nAA$ is a $P$-name in \ref{2.2}(1) 
and if in (3) we have $\Vdash \name{A_0,} \dots , \name{A_{k-1}} \in 
\nAA_2 \setminus \nAA_1$.
\end{remark*}

\begin{claim}\label{2.3}
$\leq^\kappa_{app}$ is a partial order.
\end{claim}

\proof We check (3) clause (d) of the definition. Let
$(P_1,\nAA_1) \leq^\kappa_{app} (P_2,\nAA_2)$
and $(P_2,\nAA_2) \leq^\kappa_{app} (P_3,\nAA_3)$.
Let $k < \omega$, $\name{f}_\ell$ be
$P_1$-names of injective functions from $\omega$ to $\omega$ .
Let $G \subseteq P_3$ be generic over ${\bf V}$. 
%
So let $A_\ell \in \nAA_3[G] $ for $\ell < m$.
We assume that for $\ell < m_0 \leq m$ we have that $\name{A_\ell}
\in \nAA_2$ and that  that $\{\name{A_\ell}
\such \ell < m \} \subseteq \nAA_3\setminus \nAA_2$. By the assumptions on
$P_1$ we have that $B_1 =
\left\{ n < \omega \such \bigwedge_{\ell < k }
f_\ell(n)
\not\in \bigcup \{ A_\ell \such \ell < m_0 \}\right\}$
is infinite. It belongs to ${\bf V}[G \cap P_2]$. Since we have that
$(P_2,\nAA_2) \leq_{app}^\kappa (P_3,\nAA_3)$ and 
$\{ \name{A_\ell} \such \ell \in [m_0,m) \} \subseteq
\nAA_3 \setminus \nAA_2$ and $B_1, f_0, \dots, f_{k-1} \in 
{\bf V}[G\cap P_2]$, by Definition~\ref{2.2}(3) clause (d) we are done.

\begin{claim}\label{2.4}
If $\langle (P_i, \nAA_i) \such i < \delta \rangle$ is a 
$\leq^\kappa_{app}$-increasing continuous sequence 
(continuous means that in the 
limit steps we take unions), then 
$(P,\nAA) = (\bigcup_{i < \delta}P_i, \bigcup_{i < \delta} \nAA_i)$ is an
$\leq^\kappa_{app}$-upper bound of the sequence, in particular, 
a $(\ast,\kappa)$-approximation.
\end{claim}

\proof 
The only problem is ``$(P,\nAA)$ is a $\kappa$-approximation.''

Case 1: $\cf(\delta) > \aleph_0$. Let $k < \omega$, $\name{f}_\ell$ be
$P$-names of injective functions from $\omega$ to $\omega$ .
So for some $i < \delta$ we have that 
$\langle \name{f_\ell} \such \ell < k \rangle$ is a $P_i$-name.
Let $G \subseteq P$ be generic over ${\bf V}$. In ${\bf V}[G \cap P_i]$, 
there is
some $\nAA' \subseteq \nAA$ such that
$\nAA' \in ([\nAA_i[G\cap P_i]]^{<\kappa})^{{\bf V}[G\cap P_i]}$ as 
required in ${\bf V}[G\cap P_i]$
 for $\langle \name{f_\ell}[G\cap P_i] \such \ell < k \rangle$.
We shall show that  $\nAA'$  is as required in ${\bf V}[G]$
 for $\langle \name{f_\ell}[G\cap P_i] \such \ell < k \rangle$.
So let $A_\ell \in \nAA[G] $ for $\ell < m$, w.l.o.g.\ $\name{A_\ell}
\in \nAA$, $A_\ell = \name{A_\ell}[G]$.
We assume that for $\ell < m_0 \leq m$ we have that $\name{A_\ell}
\in \nAA_i$ and that $j< \delta$ is such that $\{\name{A_\ell}
\such \ell < m \} \subseteq \nAA_j$. By the assumptions on
$P_i$ we have that $B_1 =
\left\{ n < \omega \such \bigwedge_{\ell < k }
f_\ell(n)
\not\in \bigcup \{ A_\ell \such \ell < m_0 \}\right\}$
is infinite. It belongs to ${\bf V}[G \cap P_i]$. Since we have that
$(P,\nAA_i) \leq_{app}^\kappa (P_j,\nAA_j)$ and 
$\{ \name{A_\ell} \such \ell \in [m_0,m) \} \subseteq
\nAA_j \setminus \nAA_i$ and $B_1, f_0, \dots, f_{k-1} \in 
{\bf V}[G\cap P_i]$, by Definition~\ref{2.2}(3) clause (d) we are done.

Case 2: $\cf(\delta) = \aleph_0$.
W.l.o.g.\ $\delta = \omega$. So let $k< \omega$, $p \in P$,
$p \Vdash \mbox{`` for } \ell < k, \name{f_\ell} \in {}^\omega \omega \mbox{ is injective.''}$
By renaming we may assume w.l.o.g.\ that $p \in P_0$. 
For every $m< \omega$ we find $\langle
\name{f_\ell^m} \such \ell < k \rangle$ such that 
\begin{myrules}
\item[$(\ast)_1$] $\name{f_\ell^m}$ is a $P_m$-name
for a $P/G_m$-name for an injective  function from $\omega$ to $\omega$,
\item[$(\ast)_2$] if $p \in G_m \subseteq P_m$, $G_m$ generic over ${\bf V}$ 
and
$m< \omega$, then for densely many $q \in P/G_m$ we have that
$p\Vdash_{P_m} \mbox{``}q\Vdash_{P/G_m} \bigwedge_{\ell < k}
(\name{f_\ell}\restriction m =
(\name{f_\ell^m}[G_m])\restriction m)\mbox{''}$.

\end{myrules}

So easily $p \Vdash_{P_m} \mbox{``}\name{f^m_\ell} \in 
{}^\omega A $ is injective'' where
$A$ is a countable set such that ${}^\omega A$ is the set of
all functions from $\omega$ into a set of maximal antichains
for $P/G_m$ names for functions from $\omega$ to $\omega$.
(Since we have the c.c.c.\ it is possible to make such an identification.
Also in $\nAA_m$, $\nAA'_m$, $\name{A_i}$ such an identification is
made.)
 and by the hypothesis 
on $P_m$ we have that
$p\Vdash_{P_m} \mbox{``there is } \nAA_m \in 
[\nAA_m]^{<\kappa}$ as in \ref{2.2}(1)''.
As $P_m$ is c.c.c. and because of the form of $\nAA_m$ there is $\nAA'_m$ 
a set of $<\kappa$ names from $\nAA_m$ 
such that
\begin{equation}\tag{$\ast$}
\begin{split}
& \mbox{if } \name{A_0}, \dots, \name{A_{k-1}} \in 
\nAA_m \setminus \nAA'_m \mbox{ then }\\
& p \Vdash_{P_m} \mbox{``}\left\{ 
n \such \bigwedge_{\ell < k} \name{f_\ell^m}(n)
\not\in \name{A_0} \cup \cdots \cup \name{A_{k'-1}} \right\}
\mbox{ is infinite.''}
\end{split}
\end{equation} 

So it is enough to show that $\nAA' =
\bigcup_{m < \omega } \nAA'_m$ is as required.
Let $k' < \omega$, $\name{A_0}, \dots, \name{A_{k'-1}} \in 
\nAA \setminus \nAA'$ and towards a contradiction assume that
$q \Vdash \mbox{``} \{n < \omega \such 
\bigwedge_{\ell < k} 
\name{f_\ell}(n)
\not\in \name{A_0} \cup \cdots \cup \name{A_{k'-1}} \} 
\subseteq [0,m^\ast]$.''
So for some $m$ we have that  $q \in P_m$, 
$\name{A_0}, \dots , \name{A_{k'-1}} \in \nAA_m
\setminus \nAA'_m$.
Let $q \in G_m \subseteq P_m$ be $P_m$ generic over ${\bf V}$.
In ${\bf V}[G_m]$ we have that $B' =
\{ n  \in \omega \such \bigwedge_{\ell < k}\name{f_\ell^m}[G_m](n)
\not\in \name{A_0}[G_m]\cup\dots \cup \name{A_{k'-1}}[G_m] \}
$ is infinite. So we can find $n \in B'$ such that $n > m^\ast$. 
Now there are densely many $q' \in P/G_m$ forcing $\name{f_\ell}(n)=
\name{f_\ell^m}(n)$, so w.l.o.g.\ $q \leq q' \in P/G_m$,
and we find  $p' \in G$ such that $p \leq p' \in P$ and
$p' \Vdash\mbox{``}\name{f_\ell}(n) = \name{f_\ell^m}(n)\mbox{''}.$
Contradiction.
\proofend

\begin{claim}\label{2.5}
Assume that $(P,\nAA)$ is a $\kappa$-approximation.
\begin{myrules}
\item[(1)] If 
$\Vdash \mbox{``}\name{Q} \mbox{ is Cohen or just $< \kappa$-centred }
\mbox{''}$, then $(P\ast\name{Q}, \nAA)$ is a $\kappa$-approxi\-ma\-tion,
and $(P,\nAA) \leq^\kappa_{app} (P \ast  \name{Q}, \nAA)$.
\item[(2)] If in addition $\Vdash_P \mbox{``}
\langle w_n \such n < \omega \rangle$ is a set of finite non-empty
pairwise disjoint subsets of $\omega$'', and $Q$ is Cohen forcing, and 
$\name{\eta}$ is the $P\ast \name{Q}$-name of the generic, then
$(P\ast \name{Q},\nAA \cup \{ \bigcup\{w_n \such \name{\eta}(n)=1\}\})$ is
a $\kappa$-approximation, and $\leq^\kappa_{app}$-above $(P,\nAA)$.
\end{myrules}
\end{claim}
\proof
(1) Let $G \subseteq P$ be $P$-generic over ${\bf V}$. We work in 
${\bf V}[G]$.
It is enough to prove that in $({\bf V}[G])^Q$,
$\AAA = \nAA[G]$ is a $(\kappa,\gro)$-witness.
let $Q = \bigcup_{m\in \mu} Q_m$, $Q_m$ directed, $\mu < \kappa$.  So let $\Vdash_Q
\mbox{``}\name{f_0},\dots \name{f_{k-1}} \in {}^\omega 
\omega$ are injective.''
For each $m$ we find $\langle f_\ell^m \such \ell < k \rangle$ such that 
\begin{myrules}
\item[$(\ast)_1$] $f^m_\ell \in {}^\omega \omega$,
\item[$(\ast)_2$] if $q \in Q_m$, $m < \omega$ then $q \not\Vdash_Q
\mbox{``}\neg \bigwedge_{\ell < k} \name{f_\ell}
\restriction m = f_\ell^m \restriction m\mbox{''}$.
\end{myrules}

For $\langle f_\ell^m \such \ell < k \rangle$ we have that $\AAA'_m
\in [\AAA]^{<\kappa}$ as required in Definition~\ref{2.1}(1).
 Let $\AAA'= \bigcup_{m < \mu} \AAA'_m$, it is clearly as required.

(2) We  prove clause (d) of \ref{2.2}(3).
Let $G \subseteq P$ be $P$-generic over ${\bf V}$.
So let $f_0,\dots,f_{k-1} \in V[G]$, $B \in([\omega]^\omega)^{{\bf V}[G]}$ 
and we should prove that $\{ n \in B \such 
\bigwedge_{\ell < k}
f_\ell(n) \not\in \bigcup \{ w_m \such \name{\eta}[G](n) = 1
 \}\}$ is infinite.
As $\name{\eta}$ is Cohen and the $w_n$ are pairwise disjoint and 
finite and non-empty, this follows from a density argument.
\proofend

\smallskip

An ultrafilter $D$ on $\omega$ is called Ramsey iff for
every function $f \colon  \omega \to \omega$
there is some $A \in D$ such that $f \restriction A$ is
injective or is constant.

\begin{claim}\label{2.6}
Assume that 
\begin{myrules}
\item[(a)] ${\bf V} \models \CH$,
\item[(b)] $P=\langle (P_i, \nAA_i) \such i \leq \delta \rangle$ is 
$\leq_{app}^{\aleph_1}$-increasing and continuous and $|P_i| \leq \aleph_1$,
\item[(c)] $\cf(\delta) = \aleph_1 = |\delta|$,
\item[(d)] $\delta = \sup\{i < \delta \such P_{i+1} 
= P_i \ast \mbox{\rm Cohen}, \nAA_{i+1} = \nAA_i \}$.
\item[(e)]  $G\subseteq P_\delta$ is $P_\delta$-generic over ${\bf V}$, and
in ${\bf V}[G]$ we have $\AAA = \bigcup_{i< \kappa}\nAA_i[G]$.
\end{myrules}
Then 
\begin{myrules}
\item[(1)] In ${\bf V}[G]$ there is $\bar{M}$ leisurely 
exemplifying $\AAA$.
\item[(2)] In ${\bf V}[G]$ there is a Ramsey ultrafilter $D$
such that for every $f \in {}^\omega \omega$ which is not constant on
any set in $D$ and  for all but countably [$<\kappa$] 
many $A \in \AAA$ we have
that $\{ n \such f(n) \not\in A \} \in D$.
In short we say ``$D$ is $\AAA$-Ramsey [$(\kappa,\AAA)$-Ramsey]''.
\end{myrules}
\end{claim}
\proof
(1) By renaming, w.l.o.g.\ 
$\delta= \aleph_1$.
Let $\chi\geq (2^{\aleph_0})^+$ and let $\bar{M}^0=
\langle M^0_i \such i < \omega_1 \rangle$ be increasing and
continuous and $M^0_i \prec (H(\chi), \in , <^\ast_\chi)$,
$M^0_i$ countable and $\bar{M}^0 \restriction (i+1) \in M^0_{i+1}$
and such that ${\cal P}(\omega) \subseteq \bigcup_{i<\omega_1}M_i^0$.
Let $M_i^1 = M^0_i[G]$, $\AAA_i =
\nAA_i[G]$. For any $i < \omega_1$ we shall find $j(i)
\geq i$ and $N_{j(i)}$ such that

\begin{equation}\tag{$\ast$}
\begin{split}
(\alpha) &\; M^1_{j(i)} \subseteq N_{j(i)} \subseteq M_{j(i)+1}^1,\\
(\beta) &\; N_{j(i)} \models |\AAA_{j(i)}| = \aleph_0,\\
(\gamma) &\; N_{j(i)} \in M_{j(i)+1}^1,\\
(\delta) &\; \AAA_\delta \cap N_{j(i)} = \AAA_\delta \cap M^1_{j(i)},\\
(\eps) &\; (\name{f} \in \bigcup_{i < \omega_1}
M_i^0 \; \wedge \;  \name{f}[G] \in M^1_i \cap {}^\omega \omega) 
\rightarrow \name{f}
\mbox{ is a } P_{j(i)}\mbox{-name},\\
(\zeta) &\; 
M^1_i \models |X| < \aleph_1 \Rightarrow X \subseteq M^1_{j(i)}. 
\end{split}
\end{equation}

In $M^1_i$, choose $j=j(i)$ according to the premise (d) such that
$\sup(M^1_i \cap \omega_1) < j < \omega_1$ and
$P_{j+1} = P_j \ast \mbox{Cohen}$, $\nAA_{j+1} = \nAA_j$
and such that $(\eps)$ and ($\zeta$) are true. In $M^0_{j+1}$ 
we define the forcing notion $R_j
= \{ g \such g$ is a function from some $n< \omega$ into
$\nAA_{j+1}\}$. This is a variant of Cohen forcing, and hence 
we can interpret $R_j$ as the Cohen forcing in $P_{j+1}$.
We let $\hat{g}$ be generic and set $N_j = M^1_j[\hat{g}]$. 
Now we take a club $C$ in $\omega_1$ such that
$(\forall \alpha \in C) (\forall \beta < \alpha) (j(\beta) < \alpha)$.
We let $\langle c(i)\such i < \omega_1\rangle$ be an increasing 
enumeration of $C$.
Finally we let for $i < \omega_1$,
$M_i = M^1_{c(i)}$ for limit $i$.

We have to show that in ${\bf V}[G]$,
 $\bar{M}$ $\kappa$-exemplifies ${\mathcal A}$. That is,
according to \ref{2.1}(2):
\begin{myrules}
\item[(a)] ${\mathcal A} $ is an $(\aleph_1,\gro)$-witness,
\item[(b)] $\bar{M}= \langle M_i \such i < \aleph_1 \rangle$ is 
$\prec$-increasing and continuous, and $\omega + 1 \subseteq M_0$ and 
${\mathcal P}(\omega)
\subseteq \bigcup_{i < \kappa} M_i$,
\item[(c)] $M_i \subseteq (H(\chi),\in)$ is a model of 
$\zfc^-$ and $|M_i| < \aleph_1$ and $(M_i \models |X| < \aleph_1) 
\Rightarrow X \subseteq M_i$,
\item[(d)] $\bar{M} \restriction (i+1) \in M_{i+1}$,
\item[(e)] for non-limit $i$ 
there is ${\mathcal A}_i \in M_i$ such that 
${\mathcal A} \cap M_i = {\mathcal A}_i$,
\item[(f)] if $i < \aleph_1$, $k < \omega$ and 
$f_\ell \in M_i$ is an injective function from $\omega$ to $\omega$ 
for $\ell < k$, and $k' < \omega$, $A_\ell \in 
{\mathcal A} \setminus M_i$ for $\ell < k'$, then 
$$\Bigl\{ n \such \bigwedge_{\ell < k} f_\ell(n) \not\in
A_0 \cup \cdots \cup A_{k'-1} \Bigr\} \mbox{ is infinite.}
$$
\end{myrules}
Item
(a) follows from \ref{2.4}.
The items 
(b) and (c) follow from 
 $M^0_i \prec \mbox{$(H(\chi), \in, <^\ast_\chi)$}$,
$M^0_i$ countable and $\bar{M}^0 \restriction (i+1) \in M^0_{i+1}$
and such that ${\cal P}(\omega) \subseteq \bigcup_{i<\omega_1}M_i^0$.

The items (d) and (e) are clear by our choice of $M_i$.

To show item
(f), suppose that $i < \omega_1$ and $f_\ell \in M_i$ for $\ell < k$ and
$A_\ell \in \AAA \setminus M_i$.
Then we have that $\name{f_\ell} \in V^{P_i}$ and
$A_\ell \in \AAA \setminus \AAA_i$ (the latter holds by
$(\delta)$)  and $\AAA_i = \nAA_i[G]
= \nAA_i[G_i]$ by our choice of $C$.
Hence we may use $(P_i, \nAA_i)\leq^{\aleph_1}_{app} (P_{\omega_1},\nAA)$
and get from \ref{2.2}(3)(d)
if $k < \omega$ and $\name{A_0}, \dots , \name{A_{k-1}} \in \nAA \setminus
\nAA_i$ then 
\begin{equation*}
\begin{split}
\Vdash_{P_{\omega_1}} &\mbox{`` if $B \in 
([\omega]^{\aleph_0})^{V^{P_i}}$},\\
& \name{f_\ell} \in ({}^B \omega)^{V^{P_i}} \mbox {for } \ell < k,
\mbox{ then }\\
& \left\{ n \in B \such \bigwedge_{\ell <k} \name{f_\ell}(n) 
\not\in \bigcup_{\ell <k} 
\name{A_\ell} \right\} \mbox{ is infinite'',}
\end{split}
\end{equation*}
so we get the desired property in ${\bf V}[G]$.

\nothing{
Case 1: $\cf(\delta) > \aleph_0$. Let $k < \omega$, $\name{f}_\ell$ be
$P$-names of injective functions from $\omega$ to $\omega$ .
So for some $i < \delta$ we have that 
$\langle \name{f_\ell} \such \ell < k \rangle$ is a $P_i$-name.
Let $G \subseteq P$ be generic over ${\bf V}$. In 
${\bf V}[G \cap P_i]$, there is
some $\nAA' \subseteq \nAA$ such that
$\nAA' \in ([\nAA_i[G\cap P_i]]^{<\kappa})^{{\bf V}[G\cap P_i]}$ as 
required in ${\bf V}[G\cap P_i]$
 for $\langle \name{f_\ell}[G\cap P_i] \such \ell < k \rangle$.
We shall show that  $\nAA'\subseteq M_i$  is as required in ${\bf V}[G]$
 for $\langle \name{f_\ell}[G\cap P_i] \such \ell < k \rangle$.
So let $A_\ell \in \nAA[G] $ for $\ell < m$, w.l.o.g.\ $\name{A_\ell}
\in \nAA$, $A_\ell = \name{A_\ell}[G]$.
We assume that for $\ell < m_0 \leq m$ we have that $\name{A_\ell}
\in \nAA_i$ and that $j< \delta$ is such that $\{\name{A_\ell}
\such \ell < m \} \subseteq \nAA_j$. By the assumptions on
$P_i$ we have that $B_1 =
\left\{ n < \omega \such \bigwedge_{\ell < k }
f_\ell(n)
\not\in \bigcup \{ A_\ell \such \ell < m_0 \}\right\}$
is infinite. It belongs to ${\bf V}[G \cap P_i]$. Since we have that
$(P,\nAA_i) \leq_{app}^\kappa (P_j,\nAA_j)$ and 
$\{ \name{A_\ell} \such \ell \in [m_0,m) \} \subseteq
\nAA_j \setminus \nAA_i$ and $B_1, f_0, \dots, f_{k-1} \in 
V[G\cap P_i]$, by Definition~\ref{2.2}(3) clause (d) we are done.}

(2) We work in ${\bf V}[G]$. We take
 $\langle M_i \such i < 
\omega_1 \rangle$ as in (1),
 and choose by induction on $i < \omega_1$  sets $B_i$ 
such that 
\begin{myrules}
\item[($\alpha$)] $B_i \in M_{i+1}$,
\item[($\beta$)] $j < i \Rightarrow B_i \subseteq ^\ast B_j$,
\item[($\gamma$)] if $i = j+1$ and $f \in M_j \cap {}^\omega \omega$ is
injective and $A \in \AAA \cap (M_i \setminus M_j)$, then $
B_i \subseteq^\ast \{
n \such f(n) \not\in A \} \in D$,
\item[($\delta$)] if $i$ is limit and $f \in M_i \cap {}\omega^\omega$ 
then for some
$n^\ast$ we have that $f \restriction (B_i \setminus n^\ast)$ is constant or
$f \restriction (B_i \setminus n^\ast)$ is injective.
\item[($\eps$)] $B_i$ is $<^\ast_\chi$-first of the sets
fulfilling ($\alpha$) -- ($\delta$).
\end{myrules}

Now it is easy to carry out the induction
and to show that $D$, the filter generated by
$\{B_i \such i < \omega_1 \}$ is as required.
We use property (f) of $\bar{M}$ in order to show that
requirement ($\gamma$) is no problem.
\proofend
\begin{claim}\label{2.7}
Assume that in ${\bf V}$
\begin{myrules}
\item[(a)] $\AAA$ is a $(\kappa,\gro)$-witness,
\item[(b)]  $D$ is a $(\kappa,\AAA)$-Ramsey,
\item[(c)] $Q_D = \{(w, A) \such w \in [\omega]^{<\omega} , A \in D \}$,
$(w,A) \leq (w',A')$ iff $w \subseteq w' \subseteq w \cup A$ and $A'
\subseteq A$.
\end{myrules}
Then $\Vdash_{Q_D} \mbox{``}\AAA \mbox{ is a } (\kappa,\gro) 
\mbox{-witness.''}$.
\end{claim}

\proof For $u \in [\omega]^{<\aleph_0}$ let $Q_u =
\{ (u,A) \such A \in D \}$. This is a directed subset and we have that $Q_D
= \bigcup\{Q_u \such u \in [\omega]^{<\aleph_0}\}$. So assume w.l.o.g.\
that 
$$\Vdash_{Q_D} \mbox{``}\name{f_\ell} \in 
{}^\omega \omega \mbox{ is injective for $\ell < k$''.}$$
For every $u \in [\omega]^{<\aleph_0}$ we define $f_\ell^u 
\in {}^\omega (\omega +1)$ as follows:

\begin{equation}\tag{$\otimes$} 
\begin{split}
f_\ell^u(n)= m &\mbox{ if } (\exists p \in Q_u) (p
\Vdash \name{f_\ell}(n)= m),\\
f_\ell^u(n)= \omega & \mbox{ if } (\forall m) \neg (\exists p \in Q_u)
(p \Vdash \name{f_\ell}(n)= m).
\end{split}
\end{equation}

Since $D$ is Ramsey \cite{mathias:happy} (without Ramsey
but using memory \cite{Sh:707})
we have that $Q_D$ has the pure decision property:
As usual we write $p || \varphi$ if $p\Vdash \varphi$ or $p \Vdash 
\neg \varphi$ and $q \geq_{tr} p$ iff $q \geq p$ and
$q= (w^q,A^q)$, $p=(w^p,A^p)$ and $w^q=w^p$.
\begin{equation*}\begin{split}
&
\forall p \in Q_D \exists q \geq_{tr} p
\forall u \in [\omega]^{<\aleph_0} \; \forall\ell <k \;
\forall m \in \omega \; \forall n \in (\omega+1) \\ 
&\biggl((\exists q' \geq q, q' || f^u_\ell(n)= m)
\rightarrow (\exists s \in q) (q^{[s]} || f^u_\ell(n) = m \biggr).
\end{split}
\end{equation*}
Since $Q_D$ has pure decision and $Q_u$ is directed we have that
\begin{equation}
\tag{$\ast$}
\begin{split}
&
\mbox{
for every $u \in {\omega}^{<\omega}$ for every
$m_1,m_2 < \omega$ there is some $p\in Q_u$ such that}\\
& p \Vdash \mbox{``}(\forall m < m_1) 
\min(m_2,\name{f_\ell}(m)) =\min(m_2,f_\ell^u(m)).
\end{split}
\end{equation}

For every $u \in [\omega]^{<\omega}$ and $\ell < k$ we 
can find $g^u_\ell \in {}^\omega \omega$ injective, such 
that if $\{ n \such f_\ell^u(n) < \omega \} \in D$ and $(\neg(\exists A \in D)
f^u_\ell\restriction A$ is constant) then $\{ n \such f_\ell^u(n) =
g_\ell^u(n) \} \in D$.

We call $u$  $(v,n)$-critical if
\begin{equation}\tag*{$(\ast)^u_{v,n}$}
\begin{split}
(\alpha) \; & u \in [\omega]^{<\omega},\\
(\beta) \; & \emptyset \neq v \subseteq \{0,\dots,k-1\},\\
(\gamma) \; & \ell \in v \Rightarrow f^u_\ell (n) = \omega,\\
(\delta) \; & \{m \such (\forall \ell \in v) f^{u \cup\{m\}}_\ell(n) < \omega
\} \in D ,\\
(\eps) \; & \ell < k \wedge \ell \not\in v 
\rightarrow \{ m \such f_\ell^{u \cup \{m\}}(n) = f_\ell^u(n) \} \in D.
\end{split}
\end{equation}
For $u$ $(v,n)$-critical and $\ell \in v$  
note that $\lim_D\langle f^{u\cup\{m\}}_\ell(n)
\such m < \omega \rangle = \infty$.

As $D$ is Ramsey for some $A= A_{u,v,n} \in D$ we have
if $\ell \in v$ then $\langle f_\ell^{u \cup \{m \}}(n) 
\such m \in A \rangle$ is without repetition. 

So we can find  for $\ell \in v$ injective functions
$h_\ell^{u,v,n} \in {}^\omega \omega$ 
such that $\{ m \such f_\ell^{u \cup \{m\}}(n) =
h_\ell^{u,v,n}(m) \} \in D$.

For each injective function $h \in {}^\omega \omega$ we have that 
$\AAA_h = \{ A \in \AAA \such 
\{ n \such h(n) \in A \} \in D \} $ is empty or at least
of cardinality strictly less than $\kappa$.
Let $\AAA' = \bigcup\{\AAA_h \such h = g_\ell^u$ for some $\ell < h$,
$u \in [\omega]^{<\aleph_0}$ or $h=
h_\ell^{u,v,n}$ where $u$ is $(v,n)$-critical and $\ell \in v$ and $
\emptyset \neq v \subseteq k$ $\}$.
So $\AAA' \subseteq \AAA$ is of cardinality strictly less
than $\kappa$ and it is enough to prove that
if $A_0,\dots A_{k'-1} \in \AAA \setminus \AAA'$ then
$\Vdash_Q \mbox{``}\{n \such \bigwedge_{\ell < k}
\name{f_\ell}(n) \not\in A_0 \cup \cdots \cup A_{k'-1}\} 
\mbox{ is infinite''}$.

Let $A_0,\dots,A_{k'-1}$ be given. 
Set $B^\ast= A_0 \cup \cdots \cup A_{k'-1}$.
Towards a contradiction we assume that 
$p^\ast \in Q_D$ and $n^\ast < \omega$ and
$$
p^\ast \Vdash\mbox{``}(\forall n) \left(n^\ast < n < \omega
\rightarrow \bigvee_{\ell < k} \name{f_\ell}(n) \in B^\ast\right)\mbox{''}.
$$
Let $M \prec (H(\chi),\in)$ be countable such that the following are elements
of $M$:
$p^\ast$, $D$, $\name{f_\ell}$ for $\ell < k$, $A_\ell$ for $\ell < k'$,
$\AAA'$, $\langle g^u_\ell \such u \in [\omega]^{<\aleph_0}, \ell < k
\rangle$,
$\langle h^{u,v,n}_\ell \such u \in [\omega]^{<\aleph_0}, \ell \in v,
\emptyset \neq v \subseteq k \rangle$.

Let $p^\ast=(u^\ast,A^\ast)$. Let $A^\odot \in [\omega]^\omega$ and $A^\odot 
\subseteq A^\ast$ be such that
$(\forall Y \in D \cap M)(A^\odot \subseteq^\ast Y)$ and $\min(A^\odot)
\geq\sup(u^\ast)$. It is obvious that
 $u \cup A^\odot$ is generic real for $Q_D$ over $M$, i.e.:
$\{ (u',A') \in Q_D \cap M \such u' 
\subseteq u^\ast \cup A^\odot \subseteq u' \cup A' \}$ 
is a subset of a $(Q_D)^M$-generic over $M$.

As $A_0,\dots, A_{k'-1} \in \AAA \setminus \AAA'\subseteq \AAA \setminus
\bigcup_{\ell < k} \AAA_{g_\ell^{u^\ast}}$ there is 
$n^\odot \in [n^\ast,\omega)$ such that
$\ell < k \Rightarrow g_\ell^{u^\ast}(n^\odot) \not\in B^\ast$.
Let $${\mathcal U} = 
\{ u \such u^\ast \subseteq u \subseteq u^\ast \cup A^\odot,
u \mbox{ finite, }
(\forall\ell < k)(f^u_\ell(n^\odot) < \omega 
\rightarrow f_\ell^u(n^\odot) \not\in
B^\ast \}.
$$

Now clearly $u^\ast \in {\mathcal U}$. 
Choose $u^\odot \in {\mathcal U}$ such that 
$|\{ \ell \such f^{u^\odot}_\ell(n^\odot) = \omega \}|$ is minimal.
If it is zero, we are done. So assume that is is not zero.

We choose by induction on $i< \omega$ $n_i$ such that
\begin{equation}
\tag{$\ast$}
\begin{split}
&n_i \in A^\odot,\\
&n_i < n_{i+1},\\
&\sup(u^\odot) < n_i.\\
&\ell < k \rightarrow f_\ell^{u^\odot}(n^\odot) =
f_\ell^{u^\odot \cup \{n_j \such j < i\}}(n^\odot).
\end{split}
\end{equation}

If we succeed, then $u^\odot \cup 
\{ n_i \such i < \omega \}\in M$ could have served as $A^\odot$, 
contradicting the fact that $u^\odot \cup
A^\odot$ is generic.
So for some $i$ we cannot choose $n_i$.
Let $u^\triangle = u^\odot \cup \{ n_j \such j < i \}$.
Let $v= \{ \ell < k \such \{ m \such f_\ell^{u^\triangle \cup \{m\}}(n^\odot)
\neq f_\ell^{u^\triangle}(n^\odot) \} \in D \} \subseteq 
\{0,\dots,k-1\}$.
Let $C= \{ m \such (\ell \in v \rightarrow 
 f_\ell^{u^\triangle \cup \{m\}}(n^\odot)
\neq f_\ell^{u^\triangle}(n^\odot))
\mbox{ and }
(\ell \not\in v \rightarrow 
 f_\ell^{u^\triangle \cup \{m\}}(n^\odot)
= f_\ell^{u^\triangle}(n^\odot)) \}$.
So $C \in D$ and necessarily $\ell \in v 
\wedge m \in C \Rightarrow 
 f_\ell^{u^\triangle \cup \{m\}}(n^\odot)
< f_\ell^{u^\triangle}(n^\odot)= \omega$.
So $u^\triangle$ is $(v,n^\odot)$-critical. Hence $C_1 = \{
m \such \bigwedge_{\ell \in v }h_\ell^{u^\triangle, v, n^\odot}(m) \not\in 
B^\ast \} \in D$.
Choose $n_i \in C_1 \cap C \cap M^\odot$ large enough. If $v
= \emptyset$, it can serve as $n_i$ and we have a contradiction.
Recall that $h_\ell^{u^\triangle,v,n^\odot}(n_i) = f_\ell^{u^\triangle 
\cup\{n_i\}}(n^\odot) < \infty$.
If $v \neq \emptyset$, 
then $u^\triangle\cup
\{n_i\}$ contradicts the choice of $u^\odot$, because we had
required  that
$|\{ \ell \such f^{u^\odot}_\ell(n^\odot) = \omega \}|$ is minimal.

\proofend

Later we shall use Claim \ref{1.6} in order to fulfil premise (3) of the 
following Claim~\ref{2.8}, which is together with \ref{2.3},
\ref{2.4}, \ref{2.5}, \ref{2.6} the 
justification of the single steps of our final construction of
length $\aleph_2$. Claim~\ref{2.8} serves to show 
that certain (and in the end we want to have: all) 
cofinality witnesses in intermediate $\zfc$ models
are not cofinality witnesses any more in any forcing extension. 
  
\begin{claim}\label{2.8} 
Assume that ${\bf V}$, $\langle (P_i,\nAA_i) \such i\leq \delta
\rangle$ are as in \ref{2.6}, and 
\begin{myrules}
\item[(1)] $\Vdash_{P_\delta} \mbox{``}\langle \name{K_i} \such i < \omega_1
\rangle$ is a cofinality witness and 
$\{ f \in \symom \such (\forall^\infty n) f(n) = n\}
\subseteq \name{K_0}\mbox{''.}$
\nothing{
\item[(2)] $\Vdash\mbox{``}K_i' = \name{K_i} \cap (\symom)^V
\mbox{''}$, so $\langle K_i' \such i < \omega_1 \rangle \in V$ 
is an object and
not just a name.}
\item[(2)] Let, e.g., $E_0 =
\{(n_1,n_2) \such (\exists n)(n_1,n_2 \in [n^2,(n+1)^2) \}$,
$A = \bigcup\{[(2n)^2,(2n+1)^2) \such n \in \omega \}$.
Assume that in ${\bf V}^{P_\delta}$, 
$S_{E_0,A}$ is not included in any $K_i$.
\item[(3)] $\delta = \sup \{\alpha \such \name{Q_\alpha}
\mbox{ is Cohen}, \nAA_\alpha = \nAA_{\alpha +1}\}$.
\end{myrules}

\smallskip

{\em Then} there is a $P_\delta$-name $\name{Q}$ such that
\begin{myrules}
\item[($\alpha$)] $(P_\delta,\nAA_\delta) \leq_{app}^\kappa
 (P_\delta \ast \name{Q},
\nAA_\delta)$,
\item[($\beta$)] $\Vdash_{P_\delta} \mbox{``} \name{Q} \subseteq 
\name{Q'_{E_0}}$
(where $Q'_{E_0}$ is from \ref{1.7}).
\item[($\gamma$)]
$\Vdash_{P_\delta \ast \name{Q} } 
\mbox{``}\name{g}=
\bigcup\{ \name{f} \such (p,\name{f}) \in P_\delta \ast
\name{Q} \}$ is a permutation of $\omega$ and for arbitrarily
large $i < \omega_1$,
$\langle g, \name{K_i}\rangle_{\symom}
\cap \symom^{{\bf V}[P_\delta]} \neq 
\name{K_i}$''.
\end{myrules}
\end{claim}

\proof As in \ref{2.6}, we assume w.l.o.g.\ $\delta = \omega_1$.
We can find
in ${\bf V}$, $\bar{g^\ast}=\langle \name{g_i^\ast} 
\such i < \omega_1 \rangle$ such 
that
$\Vdash_{P_{\omega_1}} \mbox{``}\name{g_i^\ast} \in 
\symom \setminus \name{K_i}$, $\name{g_i^\ast} \in S_{E_0,A}$,
$\name{g_i^\ast} \restriction (\omega \setminus A) = id$ and
$\forall n \in A$, $\name{g^\ast_i}(n) \neq n$ and $\name{g^\ast_0} \in M_0
\prec(H(\chi,\in)$, $M_0$ countable''. In ${\bf V}$ 
we now choose by induction on $i < \omega_1$  
$\name{M_i}, \name{N_i}, \name{p_i}, \alpha_i$ such that
\begin{myrules}
\item[(a)] $\langle \name{M_j} \such j \leq i \rangle$ 
is a sequence od ${\bf V}^{P_\delta}$-names  as in \ref{2.6},
\item[(b)] $\vDash_{P_\delta} \bar{Q}, \nAA, \name{\bar{g^\ast}}, 
\langle \name{K_i} \such i < \omega_1
\rangle \in \name{M_0}$,
\item[(c)] $\name{N_i}= \{ \name{\tau}_{1,n}
\such n \in \omega \}$ is a countable $P_{\alpha_i}$-name such that
$\Vdash_{P_{\alpha_i}} \mbox{``}M_i[\name{G_{P_{\alpha_i}}}] \subseteq N_i
\subseteq (H(\chi)^{{\bf V}[P_{\alpha_i}]},\in), ||N_i||=\aleph_0,
N_i \models \zfc^-\mbox{''}$,
\item[(d)] $\name{p_i} \in Q'_{E_0} $ is hereditarily countable 
and a $P_{\alpha_i}$-name of a member $Q'_{E_0}$, $
\Vdash_{P_{\alpha_i}} \langle \name{p_j} \such j \leq i 
\rangle $ is $\subseteq^\ast$-increasing and $\in \name{N_i},
\name{p_i} \in \name{N_i}$,
\item[(e)] in ${\bf V}^{P_\delta}$ we have $\name{M_i}[G_\delta] = M_i$ and
$\langle \name{N_j} \such j \leq i \rangle \in M_{i+1}$, 
$\sup(M_i \cap \omega_1)
\leq \alpha_i \in M_{i+1}$, $\name{Q_{\alpha_i}}$ is Cohen
and $\AAA_{\alpha_i}= \AAA_{\alpha_i +1}$,

\item[(f)] if $\name{I}$ is a $P_{\alpha_i}$-name of a predense subset of
$Q'_{E_0}(\langle \name{p_j} \such j < i\rangle)$, then some finite 
$J(\name{I})\subseteq \name{I}$ is predense above $\name{p_i}$ in 
$Q'_{E_0}(\langle \name{p_j} \such j \leq i\rangle)$ in the universe 
${\bf V}^{P_{\alpha_i
+1}}$.
\end{myrules}
At limit stages $i$ we take for $M_i$ the union of the former
$M_j$. Otherwise choose $M_i$ as required. Next we choose
$\alpha_i$ such that $\sup(M_i \cap \omega_1) \leq \alpha_i < \omega_1$
and $\name{Q_{\alpha_i}}$ is Cohen
and $\nAA_{\alpha_i} = \nAA_{\alpha_{i+1}}$. 
We work in ${\bf V}[P_{\alpha_i}]$.
We set $N_i^0= M_i[G_{P_{\alpha_i}}]$.
We now interpret the Cohen forcing as
$R_0 \times R_1 \times R_2$ where
$$R_0 =
\{ h \such (\exists n < \omega)
h\colon n \to {\mathcal P}(\omega)^{M_i}\}$$
ordered by inclusion.
 In $N_i^1 = N_i^0[G_{R_0}]= M_i[G_{P_{\alpha_i}}][G_{R_0}]$ we let
$$R_1 = \{ (n,q) \such n < \omega, q \in
Q'_{E_0}(\langle p_j \such j < i \rangle) \}, $$
ordered by $(n_1,q_1) \leq(n_2,q_2) \Leftrightarrow
n_1 \leq n_2 \wedge q_1 \restriction n = q_2 \restriction n \wedge
q_1 \leq q_2$. 
Since $(Q'_{E_0})^{N^1_i}$ is countable 
we have that $R_1$ is  Cohen forcing.
Let $N_i^2 = N_i^1[G_{R_0},G_{R_1}]= M_i[G_{P_{\alpha_i}}][G_{R_0}][G_{R_1}] 
$, $q_i =
\bigcup\{ q \such (n,q) \in G_{R_1}\}$.

\begin{claim*}
If $I \in {\bf V}^{P_{\alpha_i}}$ is a predense subset of 
$Q'_E(\langle p_j \such j < i \rangle)$
then for some finite $J \subseteq I$ we have:
For every $\bar{p}^\ast$ such that  $\bar{p}^\ast \restriction i =
\langle p_j \such j < j \rangle$ 
and $q_i \leq \bar{p}^\ast$ we have:
$J$ is predense above $\bar{p}^\ast$ in 
$Q'_{E_0}(\langle p_j \such j < i \rangle)$.
\end{claim*}
\proof
This is the stronger version of
\ref{1.8}(3)(b), the one starting with
``in fact \dots''. \proofend
So clearly $q_i \in (Q'_{E_0})^{{\bf V}[P_{\alpha_i +1 }]}$, 
$\bigwedge_{j<i} p_j \subseteq^\ast q_i$.

We can find in $N^2_i $ a sequence $\langle w^i_k \such k < \omega \rangle$
and $h_i^\ast$ such that
\begin{equation*}
(\ast) \left\{ \begin{array}{l}
k_1 \neq k_2 \Rightarrow w^i_{k_1} \cap w^i_{k_2} = \emptyset,\\
w^i_k \mbox{ is included in some $E_0$-equivalence class},\\
w^i_k \subseteq \omega \setminus \dom(q_i),\\
\forall n \exists m \biggl(\biggl|m/E \setminus \dom(q_i) 
\setminus \bigcup_{k \in \omega} w^i_k \biggr| > n\biggr),\\
h^\ast_i \in \symom,\\
h^\ast_i \mbox{ maps } \{n/E_0 \such n \in A\} \mbox{ onto }
\{w_k^i \such k < \omega \}\\
\mbox{more precise, $\hat{h^*_i}$ does this, where
for $b \subseteq \omega$, $\hat{h^*_i}(b) = \rge(h_i^*\restriction b)$}.
\end{array}
\right.
\end{equation*}
Let $$R_2 = \left\{ f 
\such (\exists m < \omega) \left( f \mbox{ is a permutation of }
\bigcup_{k<m} w^i_k \mbox{ mapping } w_k^i 
\mbox{ into itself}\right)\right\},$$
ordered by inclusion.
In $N^3_i = N^2_i[G_{R_2}]$ let
$f_i^\odot =
\bigcup G_{R_2}$ so $N^3_i = N^2_i[f_i^\odot]$.

So $N^3_i \in {\bf V}^{P_{\alpha_i}+1}$, and hence is a 
$P_{\alpha_i +1}$-name. As $P_{\alpha_i +1}$ has the c.c.c.,
we can assume that this name is hereditarily countable.
\nothing{We set
$q_i^0=q_i \cup f^\odot_i$ and
If $N^3_i$, $q_i^0$ are as required, 
we choose $(\name{N_i},\name{p_i})$ 
as $(N_i^3,q^0_i)$ i.e.\ $N_i^3,q_i^0$ 
are $P_{\delta}$-names so $(\name{N_i},\name{p_i})$ 
are equal to $(N^3_i,q^0_i)$.}
Now $N_i^3 \cap \omega_1 =
N^0_i \cap \omega_1 = M_i[G_{\alpha_i}] \cap \omega_1
= \delta_i < \omega_1$, hence
$N^3_i \cap \symom^{{\bf V}[P_\delta]} \subseteq K_{\delta_i}$.
Let 
$$
f^\boxdot_i = (h_i^\ast \circ g_{\delta_i}^\ast \circ (h^\ast_i)^{-1}
\restriction \bigcup_{k < \omega} w_k^i) \circ f_i^\odot.
$$
It is still generic for $R_2$ over ${\bf V}^{P_{\alpha_i}}[G_{R_0},G_{R_1}]$.
We set $N_i^4 = N^3_i[f^\boxdot_i]$,
$q'_i = q_i \cup f^\boxdot_i$.
Now $(N^4_i,q^4_i)$ are as required and choose
by taking $P_{\omega_1}$-names   $(\name{N_i},\name{p_i})$
in ${\bf V}$ for them:

Item  $(\alpha)$ of the conclusion is seen 
as follows: We have for $i < \omega_1$ that ${\bf V}^{P_{\omega_1}}
\models \mbox{``}\name{Q'_{E_0}}(\langle \name{p_j} \such j 
< i\rangle)$ is c.c.c.''. Hence we have by \ref{2.5}
that
$(P_\delta,\nAA_\delta) \leq_{app}^\kappa
 (P_\delta \ast \name{Q'_{E_0}}(\langle \name{p_j} \such j 
< i\rangle), \nAA_\delta)$, and
$(P_\delta \ast \name{Q'_{E_0}}(\langle \name{p_j} \such j 
< i\rangle),\nAA_\delta) \leq_{app}^\kappa
 (P_\delta \ast \name{Q'_{E_0}}(\langle \name{p_j} \such j 
< k\rangle), \nAA_\delta)$ for $i < k \in \omega_1$.
Since 
$ \name{Q} = \name{Q'_{E_0}}(\langle \name{p_j} \such j 
< \omega_1\rangle)=
\bigcup_{i< \omega_1}\name{Q'_{E_0}}(\langle \name{p_j} \such j < i\rangle)$
we can apply \ref{2.4}.

Item $(\beta)$ of the conclusion follows from the choice of
$\name{Q}$.

For item ($\gamma$): Fix $i$. Note that $\delta_i \geq i$.
We have in ${\bf V}^{P_{\omega_1}}$
that $f_i^{\boxdot} \in K_{\delta_i} = \name{K_{\delta_i}}[G_{\omega_1}]$.
We have that $q'_i \in (Q'_{E_0})^{{\bf V}^{P_{\alpha_i}}}$ and
$$
q'_i \Vdash_{P_{\omega_1} \ast \name{Q}}
\name{g} \restriction \bigcup_{k \in \omega}
w^i_k = \name{f_i^\boxdot} \restriction \bigcup_{k \in \omega} w^i_k 
$$
and hence
\begin{equation}\label{odot}\tag{$\odot$}
q'_i \Vdash_{P_{\omega_1} \ast \name{Q}}
g_{\delta_i}^\ast \restriction A = 
(h_i^\ast)^{-1}\circ  \name{g} \circ 
(f_i^\odot)^{-1} \circ (h^\ast_i)
\restriction A,
\end{equation}
and thus, since $ g_{\delta_i} \restriction A$ contains the same 
information as $g_{\delta_i}$ since the latter is in $S_{E_0,A}$,
the equation \ref{odot} gives a witness in 
$\langle g, \name{K_{\delta_i}}\rangle_{\symom}
\cap \symom^{{\bf V}[P_{\omega_1}]} \setminus 
\name{K_{\delta_i}}$ and hence shows the inequality claimed in $(\gamma)$.
\proofend

In order to organize the bookkeeping in our final
construction of length $\aleph_2$ we use $\diamondsuit(S_1^2)$
in order to guess the names $\langle \name{K_i} \such i < \omega_1 \rangle$
of objects that we do not want to have as cofinality witnesses.
We recall $S_1^2 = \{ \alpha \in \omega_2 \such \cf(\alpha) = 
\aleph_1 \}$. A subset of $\omega_2$ is called club
(closed and unbounded) in $\omega_2$, if it is
closed under taking suprema in the ordinals and if it is unbounded in
$\omega_2$. A subset ist called stationary, 
if its complement is not a superset of a club set.

For $E \subseteq \omega_2$ being stationary in $\omega_2$
 we have the combinatorial
principle $\diamondsuit(E)$:
There is a sequence $\langle X_\delta\such \delta \in E \rangle$ such that
for every $X \subseteq \omega_2$ the set
$\{ \delta \in E \such X_\delta = X \cap \delta \}$ is stationary in
$\omega_2$. 

For more information about this and related principles and
their relative consistency we refer the reader to \cite{Devlin,
AbrahamShelahSolovay}.

\begin{conclusion}
\label{2.9}
Assume that $2^{\aleph_0}= \aleph_1$ and that $\diamondsuit_{S^2_1}$.
Then for some forcing notion $P$ of cardinality $\aleph_2$ in 
${\bf V}^P$ we have that $\gro= \aleph_1$ and 
$\cf(\symom) = \gb = \aleph_2$.
\end{conclusion}

\proof
Let $H(\aleph_2) = \bigcup_{i < \aleph_2} B_i$, $B_i$ increasing and 
continuous, $B_{i+1} \supseteq [B_i]^{\leq \aleph_0}$ and
$\langle X_i \subseteq B_i \such i \in S_1^2 \rangle$ is a 
$\diamondsuit_{S^2_1}$-sequence.
We choose by induction on $i < \aleph_2$
$(P_i,\nAA_i,d_i)$ such that
\begin{myrules}
\item[$(\alpha)$] $(P_i,\nAA_i)$ is an $\aleph_1$-approximation,
$|P_i| \leq \aleph_1$,
\item[($\beta$)] $(P_i,\nAA_i)$ is $\leq_{app}^\kappa$-increasing and
continuous,
\item[($\gamma$)] $d_i$ is a function 
from $\nAA_i$ to $\omega_1$ (here we use that 
$\nAA_i$ is a set of $P_i$-names that are forced to be distinct),

\item[$(\delta$)] if $i< \aleph_2$ and $\langle \name{w_k} \such 
k < \omega  \rangle$ is a $P_i$-name and
$\Vdash_{P_i} \langle \name{w_k} \such 
k < \omega  \rangle $ are non-empty pairwise distinct and 
$\gamma < \omega_1$ then for some 
$j \in (i,\omega_2)$ we have that $\Vdash_{P_{j+1}}
$ for some infinite $u \subseteq \omega$ and some
$\name{A} \in \nAA_{j+1}$ we have that $\bigcup_{k \in u}
w_k \subseteq \name{A} \in \nAA_{j+1} \;\wedge \; 
d_{j+1}(\name{A}) = \gamma$,

\item[$(\eps)$] 
for arbitrarily large $i < \omega_2$ we have that $\Vdash_{P_i}
\mbox{``}Q_i = Q_{D_i}$ and $D_i$ is a Ramsey ultrafilter'',

\item[$(\zeta$)] if $i \in S_1^2$ 
and 
$P_i \subseteq B_i$, $X_i$ code of the $P_i$-name
$\langle \name{K_j} \such j < \omega_1\rangle$ and 
$\Vdash_{P_i} \mbox{``}\langle\name{K_j}\such j \in \omega_1 \rangle$ is a 
cofinality witness of $\symom^{{\bf V}[P_i]}$
 and $\{f \in \symom^{{\bf V}[P_i]}$
respects $E_0$ and $\supset id_{\omega\setminus A_0}\}$ is
not included in any $\name{K_j}$'',
then 
$\Vdash_{P_{i+1}}\mbox{`` for some } f \in \symom$ for
arbitrarily large $j < \omega_1$ we have 
$\langle \name{K_j} \cup \{f\} \rangle_{\symom} \cap
\name{(K_{j+1})^{{\bf V}_i}} \neq \name{(K_j)^{{\bf V}_i}}$''.
\end{myrules}

Can we carry out such an iteration?
We freely use the existence of limits
from Claim~\ref{2.4} and that
$\leq^*_{app}$ is a partial order \ref{2.3}. The step
$i = 0$ is trivial.
So we have to take care of successor steps.

If $i = j+1$ and $ j \not\in S_1^2$ then we can use \ref{2.5} to define 
$(P_\alpha,\nAA_\alpha)$, and taking  care of clause $(\delta)$ by 
bookkeeping.

If $i = j+1$ and $ j \in S_1^2$ and the assumption of
clause ($\zeta$) holds, we apply \ref{2.8} 
to satisfy clause $(\zeta)$, using $Q'_\zeta(\langle f_\ell
\such \ell < \omega_1 \rangle)$ from there.

If $i = j+1$ and $ j \in S_1^2$ but the assumption of
clause ($\zeta$) fails (which necessarily
occurs stationarily often), we apply \ref{2.6} and \ref{2.7}.

Having carried out the induction we let $P=
\bigcup_{\alpha < \omega_2} P_\alpha$, $\nAA = \bigcup_{\alpha< \omega_2}
\nAA_\alpha$, $d= \bigcup_{\alpha< \omega_2} d_\alpha$.
So $(P,\nAA)$ is an $(\aleph_2,\aleph_1)$-approximation.
For $\gamma \in \omega_1$ we set $\nAA^{\langle \gamma \rangle}
= \{ \name{A} \in \nAA \such d(\name{A}) = \gamma \}$.
Now clearly ${\bf V}^{P_{\aleph_2}} \models 2^{\aleph_0} = 2^{\aleph_1} =
\aleph_2$. Let $G \subseteq P$ be generic. 

We show: $\Vdash_P \gro = \aleph_1$.
For $\delta < \aleph_1$ we have that $\nAA^{\langle \delta \rangle }[G]
$ is groupwise dense by clause $(\delta)$,
and always $\gro \geq \aleph_1$. So it is
enough to show that the intersection of the
$\nAA^{\langle \delta \rangle }[G]$ is empty.
Suppose that it is not, i.e. that there is some $B 
\in [\omega]^\omega$ such that  for $\delta < \omega_1$ there is some
$A_\delta \in  \nAA^{\langle \delta \rangle }[G]$
such that for all $\delta$, $B\subseteq^\ast A_\delta$.
Now let $h \colon \omega \to B$ be an injective function.
But now we have a contradiction to
``$(P,\nAA)$ is a $(\aleph_2,\aleph_1)$-approximation (see \ref{2.3}) and 
$\nAA$ is a $(\aleph_1,\gro)$-witness (\ref{2.1}(b)).

We show that $\Vdash_{P} \gb = \aleph_2$.
This follows from clause $(\eps)$.

Finally we show that $\Vdash \cf(\symom) > \aleph_1$.
Suppose that $\langle \name{K_j}[G_{\omega_2}] \such 
j< \omega_1 \rangle$ is a cofinality witness in ${\bf V}[G_{\omega_2}]$.
Then there is a club subset $C$ in $\omega_2$ such that for
$i \in C$ we have that $\langle \name{K_j}[G_i] \such 
j< \omega_1 \rangle$ is a cofinality witness in ${\bf V}[G_i]$.
By $\diamondsuit(S^2_1)$ there is some $i \in S^2_1$ such that
$X_i$ is a code of a $P_i$ name of $\langle 
\name{K_j}[G_i] \such j < \omega_1 \rangle$.
By (the analogues of) Claims~\ref{1.4} and \ref{1.6} for $Q'_E$
 and because of $\gb=\aleph_2$ and 
because of clause ($\zeta$) we get
that the sequence $\langle \name{K_j}[G_i] \such j< \omega_1 \rangle$
does not lift to a cofinality witness in ${\bf V}[G_{\omega_2}]$ such that
for all $j < \omega_1$ we have that 
 $\name{K_j}[G_i] = \name{K_j}[G_{\omega_2}]
\cap {\bf V}[G_i]$.
Hence   $\langle \name{K_j}[G_{\omega_2}] \such 
j< \omega_1 \rangle$ was no cofinality witness in ${\bf V}[G_{\omega_2}]$.
\proofend


\begin{thebibliography}{1}

\bibitem{AbrahamShelahSolovay}
Uri Abraham, Saharon Shelah, and Robert Solovay.
\newblock Squared diamonds.
\newblock {\em Fund. Math.}, 78:165--181, 1982.

\bibitem{Devlin}
Keith Devlin.
\newblock {\em Constructibility}.
\newblock Omega Series. Springer, 1980.

\bibitem{Kunen}
Kenneth Kunen.
\newblock {\em Set Theory, An Introduction to Independence Proofs}.
\newblock North-Holland, 1980.

\bibitem{mathias:happy}
Adrian Mathias.
\newblock Happy families.
\newblock {\em Ann, Math. Logic}, 12:59--111, 1977.

\bibitem{Sh:630}
Saharon Shelah.
\newblock Non-elementary proper forcing notions.
\newblock {\em Journal of Applied Analysis}, [Sh:630], submitted.

\bibitem{Sh:h}
Saharon Shelah.
\newblock {\em {Proper and Improper Forcing}}.
\newblock Springer, 1997.

\bibitem{Sh:707}
Saharon Shelah.
\newblock Tree forcings.
\newblock {\em preprint [Sh707]}, 2000.

\bibitem{thomas:gd}
Simon Thomas.
\newblock Groupwise density and the cofinality of the infinite symmetric group.
\newblock {\em Arch. Math. Logic}, 37:483 -- 493, 1998.

\end{thebibliography}

\def\germ{\frak} \def\scr{\cal}
  \ifx\documentclass\undefinedcs\def\rm{\fam0\tenrm}\fi
  \def\defaultdefine#1#2{\expandafter\ifx\csname#1\endcsname\relax
  \expandafter\def\csname#1\endcsname{#2}\fi} \defaultdefine{Bbb}{\bf}
  \defaultdefine{frak}{\bf} \defaultdefine{mathbb}{\bf}
  \defaultdefine{beth}{BETH} \def\bbfI{{\Bbb I}} \def\mbox{\hbox}
  \def\text{\hbox} \def\om{\omega} \def\Cal#1{{\bf #1}} \def\pcf{pcf}
  \defaultdefine{cf}{cf} \defaultdefine{reals}{{\Bbb R}}
  \defaultdefine{real}{{\Bbb R}} \def\restriction{{|}} \def\club{CLUB}
  \def\w{\omega} \def\exist{\exists} \def\se{{\germ se}} \def\bb{{\bf b}}
  \def\equivalence{\equiv} \let\lt< \let\gt> \def\cite#1{[#1]}

\end{document}